\documentclass[12pt,reqno]{amsart}
\usepackage{amssymb,latexsym,amscd}

\newcommand{\dmld}{Der_{M}(log D)}
\newcommand{\dmldc}{\mathfrak{Der}_{M}(log D)}

\newcommand{\omld}{\varOmega_{M}^{1}(log D)}
\newcommand{\amld}{\mathcal{A}_{M}^{1}(log D)}

\newcommand{\dr}{\nabla^{r}}
\newcommand{\dg}{\nabla^{g}}
\newcommand{\pz}{\partial_{z}}
\newcommand{\ptz}{\partial_{\tilde{z}}}
\newcommand{\p}{\partial}

\newcommand{\pti}{\partial_{t_{i}}}
\newcommand{\ptj}{\partial_{t_{j}}}
\newcommand{\ptk}{\partial_{t_{k}}}
\newcommand{\ptl}{\partial_{t_{l}}}

\newcommand{\pqi}{\partial_{q_{i}}}
\newcommand{\pqj}{\partial_{q_{j}}}
\newcommand{\calc}{\mathcal{C}}
\newcommand{\calv}{\mathcal{V}}
\newcommand{\calu}{\mathcal{U}}
\newcommand{\calo}{\mathcal{O}}
\newcommand{\tz}{\tilde{z}}
\newcommand{\tg}{\tilde{g}}
\newcommand{\mc}{\mathbb{C}}
\newcommand{\mz}{\mathbb{Z}}
\newcommand{\mq}{\mathbb{Q}}
\newcommand{\mr}{\mathbb{R}}
\newcommand{\mch}{\mathcal{H}}
\newcommand{\mcf}{\mathcal{F}}
\newcommand{\mcu}{\mathcal{U}}

\newcommand{\mct}{\mathcal{T}}
\newcommand{\mcl}{\mathcal{L}}
\numberwithin{equation}{section}

\newtheorem{theorem}{Theorem}[section]
\newtheorem{proposition}[theorem]{Proposition}
\newtheorem{remark}[theorem]{Remark}
\newtheorem{definition}[theorem]{Definition}
\newtheorem{lemma}[theorem]{Lemma}

\begin{document}
\title[A construction of Frobenius manifolds]{A construction of Frobenius manifolds with logarithmic poles and applications}

\author{Thomas Reichelt}
\address{Institut f\"ur Mathematik; Universit\"at Mannheim, Mannheim, Germany.}
\email{thomas.reichelt@math.uni-mannheim.de}

\begin{abstract}
A construction theorem for Frobenius manifolds with logarithmic poles is established. This is a
generalization of a theorem of Hertling and Manin. As an application we prove a generalization
of the reconstruction theorem of Kontsevich and Manin for projective smooth varieties with convergent
Gromov-Witten potential. A second application is a construction of Frobenius manifolds out of a variation of polarized Hodge structures which degenerates along a normal crossing divisor when certain generation conditions are fulfilled.
\end{abstract}

\maketitle

\section*{Introduction}
\noindent In this paper we define the notion of a Frobenius manifold with logarithmic poles
along a normal crossing divisor (logarithmic Frobenius manifold for short).
Let $M$ be a complex manifold and $D$ a normal crossing divisor, then
a logarithmic Frobenius manifold is a Frobenius manifold on $M \setminus D$
such that the identity field $e$ and the Euler field $E$ are logarithmic
vector fields along $D$ and the metric $g$ and multiplication $\circ$ extend holomorphically with
respect to logarithmic vector fields and the metric is nondegenerate with respect to them.\\[5pt]

\noindent If one restricts a logarithmic Frobenius manifold to an appropriate submanifold $N$ one gets a
logarithmic Frobenius type structure whose main ingredients are a
vector bundle $K:= TM\mid_{N}$, a flat connection $\nabla^r$, a Higgs field $\calc$ which both have logarithmic poles along $D \cap N$ and a vector bundle endomorphism $\calu$ which comes from multiplication with the Euler field. If some generation conditions are satisfied by the Higgs field and the endomorphism $\calu$ the logarithmic Frobenius manifold can be uniquely
reconstructed from the logarithmic Frobenius type structure (theorem \ref{U:General}).
This is a generalization of the corresponding theorem $4.5$ of \cite{HM}.
The theorem of \cite{HM} in turn generalizes a theorem of Malgrange \cite{Mal}
which treats the case $N = \{pt\}$ and a result by Dubrovin \cite{Du}
which shows that a finite number of numbers suffices to reconstruct the whole Frobenius
germ $(M,0)$ in the case $N = \{ pt \}$, the multiplication $\circ$ being semisimple and
$\calu = E \circ $ having a simple spectrum.\\[7pt]

\noindent We give two applications of this unfolding result.\\

\noindent The first is an application to quantum cohomology. We regard the even dimensional cohomology ring
$V = H^{\ast}_{even}(X,\mc)$ as a manifold with coordinates $t_i$
with respect to a homogeneous basis $\{ T_i\}$ and metric $g_{ij} = \int_{X} T_i \cup T_j$.
The third derivatives of the Gromov-Witten potential provide us with the structure
constants of a multiplication on the tangent bundle of $V$.
As is well known these data give us a Frobenius manifold \cite{Ma}.
Let $T_1,\ldots, T_r$ be a basis of $H^2(X,\mc)$.
To construct a logarithmic Frobenius manifold we do a coordinate change $\varphi: t_i \mapsto q_i = e^{t_i}$ for $i =1, \ldots , r$ where the other coordinates remain unchanged.
The third derivatives of the Gromov-Witten potential are well defined on the image. We assume that they give holomorphic functions on a polydisc $M$ around the large radius limit point $p$ and degenerate at this point so that the corresponding multiplication becomes the cup product. The set $\bigcup\{q_i=0\}$ gives rise to a normal crossing divisor $D$. We show that the Frobenius manifold given on $M\setminus D$ has the demanded degeneration properties so that we retrieve a Frobenius manifold with logarithmic poles along $D$ on $M$.
Now let $W \subset V = H^{\ast}_{even}(X,\mc)$ be a minimal subspace which generates $H^{\ast}_{even}(X,\mc)$ with respect to the cup product. Then $N = \overline{\varphi(W)} \cap M$ is a submanifold of $M$ and carries a logarithmic Frobenius type structure. We will
use theorem \ref{U:General} to reconstruct the Frobenius manifold. In terms of Gromov-Witten invariants this means that they can be uniquely reconstructed from their values on the linear subspace $V^2 \times W^l$ for $l \geq 1$. In the case $W = H^2$ we recover the first reconstruction theorem of Kontsevich and Manin (\cite{KM}) in the case of smooth projective varieties with holomorphic Gromov-Witten potential.\\[7pt]

\noindent The second concerns with variation of polarized Hodge structures. Here we prove
that one can construct a logarithmic Frobenius type structure out of a variation of
polarized Hodge structures which degenerates along a normal crossing divisor $D$. If a certain
condition called $H^2$-generation holds this logarithmic Frobenius type structure can be unfolded
to give a logarithmic Frobenius manifold. The construction goes as follows:
With a VPHS comes the associated Gauss-Manin connection $\nabla$. One shows that the
cohomology bundle decomposes into a direct sum of subbundles using the Hodge
filtration $\mcf^{\bullet}$ and an opposite, $\nabla$-flat filtration $\calu_{\bullet}$ which is
constructed out of the limiting PMHS. Because $\nabla$ satisfies Griffiths-transversality one
can show that $\nabla$ maps $\calo(\mcf^p \cap \calu_p)$ to $\calo(\mcf^p \cap \calu_p) \oplus \calo(\mcf^{p-1} \cap \calu_{p-1})$.
Therefore one can split $\nabla$ into a connection $\nabla^r$ and a Higgs field $\calc$. This Higgs field
together with the polarization form will provide us with a logarithmic Frobenius type structure.
This section is inspired by a paper of P. Deligne \cite{De2}. Related to this section are papers of E. Cattani and J. Fernandez
\cite{CaFe} and of J. Fernandez and G. Pearlstein \cite{FePe}. In the first they prove a correspondence between
deformations of framed Frobenius modules and VPHS on $(\Delta^{\ast})^r$ with a limiting MHS of Hodge-Tate type and split over $\mathbb{R}$. In the second paper they prove in the cases of weight $k= 3,4,5$ that a deformation of framed Frobenius modules with certain generation properties gives rise to a Frobenius manifold. From the results given here this follows for any weight $k$.
\\[7pt]

\section{Unfoldings of meromorphic
connections}\label{S: unfold}  \vspace{20pt}

In this chapter we define Frobenius manifolds with
logarithmic poles (def. \ref{1:DeflogFrob}) and logarithmic Frobenius type structures (def. \ref{1:defFrobtype}).
In proposition \ref{1:FrobmfFrobstruc} it is stated that every logarithmic Frobenius manifold gives
rise to a logarithmic Frobenius type structure. To prove the construction theorem which
is the primary goal of this chapter another structure has to be established, namely
that of a $(logD-trTLEP(w))$-structure (def. \ref{1:deflogDtrtlepw}). In proposition \ref{FTSlogDoneone} we show that there is a correspondence between logarithmic Frobenius type structures and $(logD-trTLEP(w))$-structures.
The latter has the advantage that the objects coming from a logarithmic Frobenius type structure
are efficiently encoded in a vector bundle $H$ over $\mathbb{P}^1 \times M$ with a flat connection $\nabla$ and a pairing $P$. The connection is defined on $\mc^{\ast} \times M \setminus D$ having a pole of Poincar$\acute{e}$ rank $1$ along $\{ 0 \} \times M$ and a logarithmic pole along $\{ \infty \} \times M \cup \mathbb{P}^1 \setminus \{ 0 \} \times D$. In the case where the $(logD-trTLEP(w))$-structure comes from a logarithmic Frobenius manifold the pair $(H, \nabla)$ is a natural generalization of the first structure connection of a Frobenius manifold (see e.g. \cite{Du}).
These structures enable us to prove the construction theorem for Frobenius manifolds
with logarithmic poles by showing the existence of a universal unfolding of a $(logD-trTLEP(w))$-structure. \vspace{10pt} \\

\subsection{Definitions}
\subsubsection{Elementary sections}$ $\\
For the convenience of the reader we recall some definitions from
\cite{He1}.\\
Let $\Delta := \lbrace z \in \mc \mid\; \mid z \mid < 1 \rbrace$ and
$\Delta^{\ast} := \Delta \setminus \lbrace 0 \rbrace$. Fix a
holomorphic vector bundle $H \rightarrow (\Delta^{\ast})^{r} \times
\Delta^{m-r}$ of rank $\mu \geq 1$ with a flat connection $\nabla$
on $(\Delta^{\ast})^{r} \times \Delta^{m-r}$. We want to discuss
special sections in $H$ and extensions on $\Delta^{r} \times
\Delta^{m-r}$ of the sheaf $\mch$ of its
holomorphic sections.\\
Let $\mct_{i}$ be the monodromy with respect to $\{z_{i} = 0\}$. The monodromy
determines the bundle uniquely up to isomorphism. Let $\mct_{i} =
\mct_{i,s} \cdot \mct_{i,u}$ be the decomposition into semisimple
and unipotent parts, $N^{(i)}:= \log(\mct_{i,u})$ the nilpotent part.
A universal covering is
\begin{align}
e: \mathbb{H}^{r} \times \Delta^{m-r} &\longrightarrow
(\Delta^{\ast})^{r} \times \Delta^{m-r}, \notag \\
(\xi_{1}, \ldots , \xi_{r},z_{r+1}, \ldots, z_{m}) &\mapsto (e^{2
\pi i \xi_{1}}, \ldots , e^{2 \pi i \xi_{r}}, z_{r+1}, \ldots ,
z_{m}). \notag
\end{align}
$ $\\
The space of multivalued flat sections $H^{\infty}$ is defined as
\[
H^{\infty} = \lbrace pr \circ A : \mathbb{H}^{r} \times \Delta^{m-r}
\rightarrow H \mid A \text{ is a global flat section of}\; e^{\ast}H
\rbrace.
\]
$ $\\
All monodromies act on it. The indices of the simultaneous
generalized eigenspace decomposition $H^{\infty} =
\bigoplus_{\lambda} H_{\lambda}^{\infty}$ can be considered as
tuples $\underline{\lambda} = (\lambda_{1}, \ldots, \lambda_{r})$ of
eigenvalues
for the monodromies $\mct_{1}, \ldots, \mct_{r}$.\\
Now choose $A \in H_{\underline{\lambda}}^{\infty}$ and
$\underline{\alpha} = (\alpha^{(1)}, \ldots , \alpha^{(r)})$ with $e^{-2
\pi i \alpha^{(j)}} = \lambda_{j}$. The map
\begin{align}
\mathbb{H}^{r} \times \Delta^{m-r} &\longrightarrow H, \notag \\
\xi =(\xi_{1}, \ldots , \xi_{r},z_{r+1}, \ldots, z_{m}) &\mapsto
\prod_{j=1}^{r} \left( \exp(2 \pi i \, \alpha^{(j)} \xi_{j})
\exp(-\xi_{j} N^{(j)})\right) A(\xi) \notag
\end{align}
is invariant with respect to any shift  $\xi_{j} \mapsto \xi_{j}+1$
and therefore  induces a holomorphic section
\begin{align}
es(A, \underline{\alpha}) &: (\Delta^{\ast})^{r} \times \Delta^{m-r}
\rightarrow H, \notag \\
z= (z_{1}, \ldots, z_{m}) &\mapsto \prod_{j=1}^{r} \left( \exp(2 \pi
i \, \alpha^{(j)} \xi_{j}) \exp(-\xi_{j} N^{(j)})\right) A(\xi) \notag
\end{align}
of the bundle $H$ (with $z = e(\xi)$). It is called an elementary
section and is usually denoted informally as
\[
z \mapsto \prod_{j=1}^{r} z_{j}^{\alpha^{(j)}-\frac{N^{(j)}}{2 \pi i}}A\; .
\]
It is nowhere vanishing if $A \neq 0$ because the twist with
$\prod_{j=1}^{r} z_{j}^{\alpha^{(j)}- \frac{N^{(j)}}{2 \pi i}}$ is invertible.\\[5pt]

\noindent We now want to investigate extensions of $\mch$ by
elementary
sections.\\
A $\mz^{r}$-filtration $(P_{\underline{l}})_{\underline{l} \in
\mz^{r}}$ consists of subspaces $P_{\underline{l}} \subset
H^{\infty}$ which are invariant with respect to all monodromies
$\mct_{1}, \ldots, \mct_{r}$ and which satisfy for a suitable $m >
0$
\begin{align}
&P_{\underline{l}} = 0 \quad \text{if $l_{j} \leq -m$ for some $j$},
\notag \\
&P_{\underline{l}} = H^{\infty} \quad \text{if $l_{j} \geq m$ for
all $j$}, \notag \\
&P_{\underline{l}} \subset P_{\underline{l'}} \quad \text{if $l_{j}
\leq l'_{j}$ for all $j$}. \notag
\end{align}
Then $P_{\underline{l}} = \bigoplus_{\underline{\lambda}}
P_{\underline{l}, \underline{\lambda}}$. A $\mz^{r}$-filtration
$(P_{\underline{l}})$ induces $r$ increasing exhaustive monodromy
invariant filtrations $F_{\bullet}^{j}$ on $H^{\infty}$ by
\[
F_{p}^{j} := \bigcup_{l_{j} = p} P_{\underline{l}}\;.
\]
For $\lambda_{j}$ an eigenvalue of $\mct_{j}$ let $\alpha_{j}$ be
defined as $e^{-2 \pi i \alpha_{j}} = \lambda_{j}$ and $-1 < Re\,
\alpha_{j} \leq 0$. Set $M := \Delta^{r} \times \Delta^{m-r}$,
$M^{\ast} := (\Delta^{\ast})^{r} \times \Delta^{m-r}$ and $D:= M \setminus M^{\ast}$.
\vspace{10pt}
\begin{proposition}\label{P: es}
Fix a $\mz^{r}$-filtration $(P_{\underline{l}})_{\underline{l} \in
\mz^{r}}$.
\begin{enumerate}
\item For any $\underline{l} \in (\lbrack -m,m \rbrack \cap
\mz)^{r}$ and any $\underline{\lambda}$, choose generators
$A^{i}_{\underline{l}, \underline{\lambda}}$ of $P_{\underline{l},
\underline{\lambda}}$. The sheaf
\[
\mcl := \sum_{\text{these}\, \underline{l}}
\sum_{\underline{\lambda}} \sum_{i} \calo_{M} \cdot
es(A^{i}_{\underline{l},\underline{\lambda}}, \underline{\alpha} +
\underline{l})
\]
is an $\calo_{M}$-coherent extension of $\mch$ with logarithmic pole
along $D$.
\item The formula above yields a one-one correspondence between
$\mz^{r}$-filtrations and $\calo_{M}$-coherent extensions of $\mch$
with logarithmic pole along $D$.
\item The sheaf $\mcl$ is a locally free $\calo_{M}$-module if and
only if
\[
P_{\underline{l}} = \bigcap_{j=1}^{r} F_{l_{j}}^{j} \quad \text{for
all}\; \underline{l}
\]
holds and the filtrations $F_{\bullet}^{j}$ have a common splitting.
\end{enumerate}
\end{proposition}
\begin{proof}
See \cite{EV} Appendix C and \cite{He1} chapter 8.3 .
\end{proof}
\vspace{10pt}
\begin{remark}
A special case of proposition \ref{P: es} is the case when the
filtrations on the simultaneous generalized eigenspaces
$H_{\lambda}^{\infty}$ are all of type $ 0 =
F_{p(j,\lambda)-1}^{j} \subset F_{p(j,\lambda)}^{j}= H^{\infty}_{\lambda}$
for some numbers $p(j,\lambda) \in \mz$.\\
The $\calo_{M}$-locally free extension of $\mch$ to $M$ for the case
when the real parts of all eigenvalues of all residue endomorphisms
are contained in $\lbrack 0,1 )$ was called the canonical extension
by Deligne.
\end{remark}

\subsubsection{\textbf{(Logarithmic) Frobenius manifolds}}
\normalsize $ $\vspace{10 pt}  \\ Before we can delve into
logarithmic Frobenius manifolds, we state the definition of a
Frobenius manifold.\\
Let $M$ be a complex manifold and $TM$ its (holomorphic) tangent bundle.
In this context a metric $g$ on $TM$ is a symmetric, $\mc$-bilinear, nondegenerate, holomorphic pairing on the fibers of the holomorphic bundle $TM$ and a multiplication $\circ$ is a $\mc$-bilinear, holomorphic, commutative and associative multiplication on the fibers of $TM$.
\begin{definition}\cite{Du}\label{1:DefFrob}
A Frobenius manifold is a tuple $(M, \circ, e ,E, g)$ where $M$ is a
manifold of dimension $n \geq 1$ with metric $g$ and multiplication
$ \circ$ on the tangent bundle, $e$ is a global unit field and $E$
is another global vector field, subject to the following conditions:
\begin{enumerate}
\item the metric is multiplication invariant, $g(X \circ Y, Z) =
g(X, Y \circ Z)$,
\item (potentiality) the (3,1)-tensor $\nabla \circ$ is symmetric (
here $\nabla$ is the Levi-Civita connection of the metric),
\item the connection $\nabla$ is flat,
\item the unit field $e$ is flat, $\nabla e = 0$,
\item the Euler field $E$ satisfies $Lie_{E}(\circ) = 1 \cdot \circ$
and $Lie_{E}(g) = (2-d) \cdot g$ for some $d \in \mc$.
\end{enumerate}
\end{definition}
\noindent Now we can give  the definition of a logarithmic Frobenius
manifold.
\begin{definition}\label{1:DeflogFrob}
Let $M$ be a complex manifold and $D$ a normal crossing divisor. Let
$M \setminus D$ be a Frobenius manifold with
\begin{align}
e&, E \in \Gamma(M,\dmldc), \notag \\
& \circ \hspace{8pt} \in \Gamma(M,\dmldc \otimes ( \amld
\otimes_{S} \amld)), \notag \\
& g \hspace{8pt} \in \Gamma(M,( \amld \otimes_{S} \amld)) \notag
\end{align}
where $g$ is nondegenerate on $\dmld$.\\
Then M is a Frobenius manifold with logarithmic poles along D
(logarithmic Frobenius manifold for short).
\end{definition}

\begin{proposition}\label{1:LClp}
The induced Levi-Civita connection on $M \setminus D$ with respect
to $g$ has a logarithmic pole along $D$.
\end{proposition}
\begin{proof}
Just compute the Christoffel symbols in a neighborhood $U$ of the divisor with
local coordinates $(t_{1}, \ldots , t_{n})$ such that $D \cap U = \bigcup_{i=1}^{r} \{ t_i = 0 \}$.
\end{proof}
\noindent Define
\[
\calu := E \circ : \dmldc \rightarrow \dmldc
\]
and
\begin{align}
\calv : \; &\dmldc \rightarrow \dmldc \notag \\
X &\mapsto \nabla_{X}E - \frac{2-d}{2}X\; . \notag
\end{align}
\vspace{15pt}

\subsection{Frobenius type structures with logarithmic poles} $ $\vspace{15pt}\\
In this section we define Frobenius type structures with logarithmic
poles. The advantage of these structures is, that restricted to any
submanifold they remain Frobenius type structures with logarithmic
poles. Frobenius manifolds do not have this property in general because the multiplication $\circ$ is defined on the tangent bundle $TM$ of $M$.
\vspace{10pt}
\begin{definition}\label{1:defFrobtype}
 A structure of Frobenius type with logarithmic poles along $D$ is a
 tuple $(K \rightarrow (M,0) , D, \dr, \calc, \calu, \calv , g)$
 with
\begin{align}
& (M,0) \; \text{is a germ of a complex manifold}, \notag \\
& D\; \text{is a normal crossing divisor with}\; 0 \in D \;
\text{locally given by}\; D = \{z_{1}\cdots z_{r} \}, \notag \\
& K \longrightarrow (M,0)\; \text{germ of a holomorphic vector bundle}, \notag \\
& \dr \; \text{flat connection with logarithmic pole along D}, \notag \\
& \calc : \calo (K) \longrightarrow \calo (K) \otimes \amld
\quad \text{logarithmic Higgs field}, \notag \\
& \calu : \calo (K) \longrightarrow \calo (K) \;  \text{with} \; [\,
\calu , \calc_{X}] = 0 \quad \forall X \in \dmldc, \notag \\
& \calv : \calo(K) \longrightarrow \calo (K) \quad \dr \text{-
flat},
\notag \\
& g : \calo (K) \times \calo (K) \longrightarrow \calo_{M,0} \;
\text{symmetric, non-degenerate,}\;  \dr \text{- flat}, \notag
\end{align}
such that these data satisfy
\begin{align}
& \dr_{X} (\calc_{Y}) - \dr_{Y}(\calc_{X}) - \calc_{[X,Y]} = 0, \label{D :FTS1} \\
& \dr ( \calu ) - [\calc , \calv] + \calc = 0, \label{D :FTS2} \\
& g(\calc_{X} a,b) = g(a, \calc_{X} b), \label{1:gC}\\
& g(\calu a,b) = g(a, \calu b), \label{1:gU}\\
& g(\calv a,b) = -g(a, \calv b) \label{1:gV}.
\end{align}
\end{definition}
\vspace{15pt}
\newpage
\begin{proposition}\label{1:FrobmfFrobstruc}
Every Frobenius manifold with logarithmic poles gives rise to a Frobenius type
structure with logarithmic poles.
\end{proposition}

\begin{proof}
The proof is essentially the same as in \cite{He2} Lemma 5.11, so we just define the objects:
\begin{align}
& K := \dmld , \quad \dr := \nabla^{g}, \notag \\
& \calu := E \circ :\; \dmldc \longrightarrow  \dmldc, \notag \\
& \calv :\; \dmldc \longrightarrow \dmldc\; := \quad X \mapsto
\nabla_{X} E
- \frac{2-d}{2} X, \notag \\
& \calc : \dmldc \longrightarrow \dmldc \otimes \amld\, , \quad
\calc_{X} Y := - X \circ Y. \notag
\end{align}
\end{proof}

\subsection{Definition of (logD-trTLEP(w))-structures} $ $\normalsize\vspace{15pt}
\\
Let $(M,0)$ be a germ of a complex manifold and $D$ be a normal crossing divisor
and let $z$ be a coordinate on $\mathbb{P}^1$.

\begin{definition}\label{1:deflogDtrtlepw}
Fix $w \in \mathbb{Z}$. A $(logD - trTLEP(w))$-structure is a tuple
$((M,0),$ $D,H, \nabla, P)$ with $H \hspace{-3pt} \rightarrow \hspace{-1pt} \mathbb{P}^{1} \times (M,0)$ a
trivial holomorphic vector bundle with a flat connection on $H_{\mid
\mathbb{C}^{*} \times (M \setminus D ,0)}$ with
\[
\nabla : \calo(H)_{(0,0)} \longrightarrow \frac{1}{z} \Omega^{1}_{\mc \times M,0}(log Z_0) \otimes \calo (H)_{(0,0)}.
\]
where $(Z_0,(0,0)) := (\{ 0 \} \times M \cup \mc \times D,(0,0) ) \subset (\mc \times M, (0,0) )$
and $\nabla$ has a logarithmic pole along $\{ \infty\} \times (M,0) \cup (\mathbb{P}^1\setminus\{0\}) \times (D,0) \subset (\mathbb{P}^1\setminus\{0\}) \times (M,0) $.
$P$ is a
$(-1)^{w}$-symmetric, nondegenerate, $\nabla$-flat pairing
\[
P : H_{(z,t)} \times H_{(-z,t)} \longrightarrow \mathbb{C} \quad
\text{for}\; (z,t) \in \mathbb{C}^{*} \times (M \setminus D,0).
\]
on a representative of $H$ such that the pairing extends to a
nondegenerate,\\ $z$-sesquilinear pairing
\[
P : \calo(H) \times \calo(H) \longrightarrow z^{w} \calo_{\mathbb{P}^1
\times (M,0)}.
\]
\end{definition}

\vspace{10pt}
\noindent If $((M,0),D,$ $H, \nabla, P)$ is a $(logD-trTLEP(w))$-structure and $\varphi: (M',D',0) \rightarrow (M,D,0)$
a holomorphic map of germs of manifolds then one can pullback $H$, $\nabla$ and $P$ with $id \times \varphi : (\mc,0) \times (M',0) \rightarrow (\mc,0) \times (M,0)$.
One easily sees that $\varphi^{\ast}(H,\nabla)$ and $\varphi^{\ast}(P)$ gives a $(logD'-trTLEP(w))$-structure on $(M',D',0)$.
\begin{definition}
Fix a $(logD-trTLEP(w))$-structure $((M,0),D,H,\nabla,P)$.\\
(a) An unfolding of it is a $(logD-trTLEP(w))$-structure $((M \times
\mc^{l},0),D \times \mc^{l}, \tilde{H},$ $\tilde{\nabla}, \tilde{P})$
together with a fixed isomorphism
\[
i:((M,0),D,H,\nabla,P) \rightarrow ((M\times \mc^{l},0), \tilde{D},
\tilde{H}, \tilde{\nabla}, \tilde{P})\mid_{(M \times \{ 0 \} ,0)}
\]
where we denote $D \times \mc^{l}$ as $\tilde{D}$.\\
(b) One unfolding $((M \times \mc^{l},0),D \times \mc^{l},
\tilde{H}, \tilde{\nabla}, \tilde{P},i)$ induces as second unfolding
$((M \times \mc^{l'},0),\tilde{D}',$ $\tilde{H}', \tilde{\nabla}',
\tilde{P}',i')$ if there is an isomorphism $j$ from the second
unfolding to the pullback of the first unfolding by a map
\[
\varphi: ((M \times \mc^{l'},0), \tilde{D}') \rightarrow ((M \times
\mc^{l},0), \tilde{D}),
\]
which is the identity on $(M \times \{ 0 \},0)$, and if
\[
i = j\mid_{(M \times \{ 0 \},0)} \circ i'.
\]
(c) An unfolding is universal if it induces any unfolding via a
unique map $\varphi$.
\end{definition}

\subsection{A Correspondence} $ $ \normalsize\\
\begin{proposition}\label{FTSlogDoneone}
There is a 1-1 correspondence between
$(logD-trTLEP(w))$-structures and Frobenius type structures with
logarithmic pole.
\end{proposition}
\begin{proof}$ $\\
\textbf{Logarithmic Frobenius type structure $\rightarrow$
(logD-trTLEP(w))-structure}.\\ $ $ \\
Let $(K \rightarrow (M,0),D, \dr ,\calc, \calu, \calv, g)$ be a
Frobenius type structure with logarithmic pole. Set $H := \pi^{*}K$
and let $\varphi_{z}: H_{(z,t)} \rightarrow K_{t}$ be the canonical
projection. Extend $\dr, \calc, \calu, \calv, g$ canonically on H.
Define
\[
\nabla :=  \dr + \frac{1}{z}\calc + (\frac{1}{z} \calu - \calv +
\frac{w}{2} id) \frac{dz}{z}
\]
and
\begin{align}
P : H_{(z,t)} \times H_{(-z,t)} &\rightarrow \mathbb{C} \quad
\text{for } (z,t)
\in (\mathbb{C}^{*},0) \times (M,0), \notag \\
(a,b) &\mapsto z^{w}g(\varphi_{z}a,\varphi_{-z}b). \notag
\end{align}
We have to show that $(H,D, \nabla,P)$ is a
$(logD-trTLEP(w))$-structure.\\
$H = \pi^{*}K$ is trivial because $K$ is trivial. $\nabla^{2}= 0$
results from a quick calculation. $P$ is $(-1)^{w}$-symmetric and
nondegenerate because $g$ is symmetric and nondegenerate.\\
Let $a,b \in \calo(H)$. $P$ is $\nabla$-flat: \footnote{for better
readability we suppress $\varphi_{\pm z}$}
\begin{align}
dP(a,b) &= w z^{w-1} \cdot g(a,b) + z^{w}
\cdot dg(a,b), \notag \\
P(\nabla a, b) &= z^{w}g((\dr + \frac{1}{z} \calc +(\frac{1}{z}
\calu - \calv + \frac{w}{2}
id)\frac{dz}{z})a,b), \notag \\
P(a, \nabla b) &= z^{w}g(a,(\dr - \frac{1}{z} \calc + (-\frac{1}{z}
\calu - \calv +\frac{w}{2} id)\frac{dz}{z})b) \notag
\end{align}
\begin{align}
\Longrightarrow \nabla(P)(a,b)= &dP(a,b) - P(\nabla a, b) - P(a, \nabla b) = 0 \; .
\end{align}
Because of the definition of $\nabla$ above one sees easily that it has the
correct pole order and that therefore $(H,D, \nabla,P)$ is a $(logD-trTLEP(w))$-structure.\\[5pt]

\noindent \textbf{(logD-trTLEP(w))-structure $\rightarrow$ Logarithmic
Frobenius type structure.}\\ $ $\\
Let a $(logD-trTLEP(w))$-structure $((M,0),D,H,\nabla,P)$ be given.\\
Define $K:= H_{\mid \{0\} \times (M,0)}$. Observe that $K$ and
$H_{\mid \{\infty \} \times (M,0)}$ are canonically isomorphic.\\
The residual connection $\nabla^{res}$ on $H_{\mid \{\infty \}
\times (M,0)}$ is given by
\[
\nabla^{res}_{X} [a] := [ \nabla_{X} a ] \quad \text{with} \quad X
\in \dmldc,\; a \in \calo(H).
\]
The residue endomorphism $\calv^{res}$ on $H_{\mid \{\infty \}
\times (M,0)}$ is given by
\[
\calv^{res}[a] := [\nabla_{\tz \ptz} a] \quad a \in \calo(H).
\]
Because of the canonical isomorphism these structures can be shifted
to $K$. Let $\dr, \calv$ be the shift of $\nabla^{res},
\calv^{res}+\frac{w}{2}id$ respectively. Furthermore we set
\begin{align}
&\calc := [z \nabla] : \calo(K) \longrightarrow \calo (K) \otimes
\amld, \notag \\
&\calc_{X} [a] := [z \nabla_{X} a], \notag \\
& \calu := [z \nabla_{z \pz}] : \calo(K) \longrightarrow \calo(K),
\notag \\
& g([a],[b]) := z^{-w} P(a,b)\; \text{mod}\; z \calo_{\mathbb{C}
\times M,0}. \notag
\end{align}
$\dr$ is flat because $\nabla$ is flat and has a logarithmic pole
along $D$ because $\nabla$ has a logarithmic pole along $\hat{D}$.
\\
It follows from easy calculations that $\calc$ is a logarithmic Higgs field, $[\calu,\calc_X] =0$ and $\nabla^r(\calv) = 0$.\\

\noindent $g$ is $\calo_{M,0}$-bilinear, symmetric and nondegenerate because
$P$ is $\calo_{\mathbb{C} \times M}$ $z$-sesqui-linear,
$(-1)^{w}$-symmetric and nondegenerate. The conditions $(\ref{1:gC}), (\ref{1:gU}), (\ref{1:gV})$ and $\nabla^r(g)=0$ follow from $\nabla(P)=0$.\\
Lift $\dr, \calc, \calu, \calv$ canonically to $H$. Consider
\begin{align}
\Delta := \nabla -(\dr + \frac{1}{z} \calc +(\frac{1}{z} \calu -
\calv + \frac{w}{2} id)\frac{dz}{z}).
\end{align}
$\Delta$ is $\calo_{\mathbb{P}^{1} \times M}$-linear. $\Delta_{X}$ for $X \in \pi^{-1}\dmldc$ maps sections in
$\pi^{-1} \pi_{*} \calo(H)$ to sections with no pole along $\{ 0 \}
\times (M,0)$ because of the definition of $\calc$ and which vanish along $\{ \infty \} \times (M,0)$ because of the definition of $\nabla^r$.
Therefore $ \Delta_X = 0$. The same holds for $\nabla_{z \p_z}$, using the definitions of $\calu$ and $\calv$. This shows $\Delta = 0$.\\[5pt]
\noindent We show condition (\ref{D :FTS1}). Let $a \in
\pi^{-1} \pi_{\ast} \calo(H)$. It holds $ \nabla_{X} \nabla_{Y} -
\nabla_{Y} \nabla_{X} - \nabla_{[X,Y]} =0 $ with $ \nabla_{X} =
\dr_{X} + \frac{1}{z} \calc_{X}:$
\begin{align}
&\nabla_{X} \nabla_{Y} a = \dr_{X} \dr_{Y} a + \dr_{X} \frac{1}{z}
\calc_{Y} a+ \frac{1}{z} \calc_{X} \dr_{Y} a + \frac{1}{z^{2}}
\calc_{X} \calc_{Y} a ,\notag \\
&\nabla_{Y} \nabla_{X} a = \dr_{Y} \dr_{X} a + \dr_{Y} \frac{1}{z}
\calc_{X} a+ \frac{1}{z} \calc_{Y} \dr_{X} a + \frac{1}{z^{2}}
\calc_{Y} \calc_{X} a , \notag \\
&\nabla_{[X,Y]}a = \dr_{[X,Y]}a + \frac{1}{z} \calc_{[X,Y]}a, \notag
\\
&\Longrightarrow \frac{1}{z} (\dr_{X}(\calc_{Y}) -
\dr_{Y}(\calc_{X}) - \calc_{[X,Y]})a = 0 . \notag
\end{align}
Condition \ref{D :FTS2} follows from a similar calculation.
\end{proof}

\subsection{The isomorphism case} $ $
\\[10pt]
In this section we prove a result which shows that if the Higgs
field of a logarithmic Frobenius type structure provides an
isomorphism between the fiber over 0 and the fiber over 0 of the
tangent bundle and if a certain eigenvector condition is satisfied,
then we can construct a logarithmic Frobenius manifold out of the
logarithmic Frobenius type structure. This is a cornerstone for the
general result which we prove in the next section.
\begin{proposition}\label{1:isocase}
Let $(K \rightarrow (M,0), D , \dr, \calc, \calu,\calv, g)$ be a
Frobenius type structure with logarithmic pole along D. Let $\xi \in
\calo(K)$ be a section which is flat on $M \setminus D$ such that
\begin{align}
\calc : \dmld_{0} &\longrightarrow K_{0}\, , \notag \\
X &\mapsto - \calc_{X} \xi\hspace{-2pt}\mid_{0} \quad \text{is an isomorphism} \notag \\
\text{and}\;\; \calv \xi &= \frac{d}{2} \xi \; .\notag
\end{align}
Then the Frobenius type structure corresponds to a germ of a
Frobenius manifold with logarithmic pole along D ((M,0), D, $\dg$,
$\circ$, e, E, $\tilde{g}$).
\end{proposition}
Notice that, unlike in the corresponding result of \cite{HM} where the vector $\xi\hspace{-3pt}\mid_0$ could be
extended to a flat section, we have to demand
the existence of such a section $\xi$ which is flat on $M\setminus D$. The section $\xi$ corresponds
to the global unit field $e$ under the isomorphism given above.
\begin{proof}
We have an isomorphism
\[
v: \dmldc \longrightarrow \calo(K)\, ,
\]
so we can pullback the data of the Frobenius type structure
\[
\dg := v^{*}(\dr),\; \tg := v^{*}(g),\; e:= v^{-1}(\xi),\; E:=
v^{-1}(\calu \xi).
\]
$\dg$ is flat with $\dg \tg =0$ and $\dg e = 0$ and has a
logarithmic pole along D. Torsion freeness follows from condition $(\ref{D :FTS1})$.\\[5pt]
\noindent Define the multiplication on $\dmld$
\[
v(X \circ Y) := -\calc_{X}v(Y) = \calc_{X} \calc_{Y} \xi\; .
\]
The potentiality follows from easily:
\begin{align}
0 &= \dr_{X}(\calc_{Y}) - \dr_{Y}(\calc_{X}) - \calc_{[X,Y]}
\notag \\
\Leftrightarrow & -\dg_{X}(Y \circ) + \dg_{Y}(X \circ) + [X,Y] \circ
= 0
\notag \\
\Rightarrow & -\dg_{X}(Y \circ Z) + Y \circ \dg_{X}(Z) + \dg_{Y}(X
\circ Z) - X \circ \dg_{Y}Z + [X,Y] \circ Z = 0 \notag
\end{align}
but this is equivalent to $\dg(\circ)$ symmetric (see \cite{He1}).\\
We denote the shift of $\calv, \calu, \calc$ on
$\dmld$ with the same letters. We want to show
\[
\calv = \dg_{\bullet} E - \frac{2-d}{2}id = \dg_{E} - Lie_{E} -
\frac{2-d}{2} id\; .
\]
We need the following identities for the next calculation:
\begin{align}
\nabla_{z \pz} v(Y) &= (\frac{1}{z} \calu - \calv +\frac{w}{2}
id)v(Y) \notag \\
&= \frac{1}{z} v(E \circ Y) - v((\calv - \frac{w}{2} id)Y), \notag \\
\nabla_{X} v(e) &= v(\dg_{X}e) -\frac{1}{z} v(X) = -\frac{1}{z} v(X)
. \notag
\end{align}
$ $\\[20pt]
Therefore we have
\begin{align}
&v(\dg_{X}E) = \dr_{X}v(E) \notag \\
=\; &\nabla_{X}v(E) - \frac{1}{z} \calc_{X}
v(E) = \nabla_{X} \calu v(e) + \frac{1}{z}\calu v(X) \notag \\
= \;&\nabla_{X} z(\nabla_{z \pz} v(e) + v(\frac{d-w}{2}e)) +
(\nabla_{z \pz} v(X) + v((\calv- \frac{w}{2}id)X)) \notag \\
= \;&- \nabla_{z \pz} v(X) + v(X) - v(\frac{d-w}{2}X)+ \nabla_{z \pz}
v(X) + v((\calv - \frac{w}{2} id)X) \notag \\
= \;&v((\calv + \frac{2-d}{2} id)X) \notag
& \notag \\
&\hspace{-13pt}\Rightarrow  Lie_{E} = \dg_{E} - \calv -\frac{2-d}{2} id . \notag
\end{align}
$ $\\[10pt]
Now we can compute $Lie_{E}(g)$ and $Lie_{E}(\circ)$:
\begin{align}
& Lie_{E}(g) (X,Y) \notag \\
=\; & E g(X,Y) - g(Lie_{E}X,Y) - g(X,Lie_{E} Y)
\notag
\\
=\; & Eg(X,Y) - g(\dg_{E}X,Y) - g(X,\dg_{E}Y) + g(\calv X,Y) + g(X,
\calv Y) \notag \\
 &+ g(\frac{2-d}{2}X,Y) + g(X, \frac{2-d}{2}Y) \notag \\
=\; & (2-d) g(X,Y) . \notag
\end{align}
$ $\\[10pt]
For the next calculation we are using $\dg_{X}(\calc_{Y}) -
\dg_{Y}(\calc_{X}) - \calc_{[X,Y]} = 0$,\\ $\nabla(\calu) - [\calc,
\calv] + \calc = 0$ and $\calu = - \calc_{E}$:
\begin{align}
&\;\calc_{Lie_{E}(\circ)(X,Y) - X \circ Y} \notag \\
=& \; \calc_{[E, X \circ Y]}
- \calc_{X \circ [E,Y]} - \calc_{[E,X] \circ Y} - \calc_{X \circ Y}
\notag \\
=& \;\nabla^{g}_{E}(\calc_{X \circ Y}) - \nabla^{g}_{X \circ Y}(\calc_{E}) +
\calc_{X} \calc_{[E,Y]} + \calc_{[E,X]}\calc_{Y} +
\calc_{X}\calc_{Y} \notag \\
=& - \nabla^{g}_{E}(\calc_{X} \calc_{Y}) - \nabla^{g}_{X \circ Y}(\calc_{E})
+ \calc_{X}(\nabla^{g}_{E} (\calc_{Y}) - \nabla^{g}_{Y} (\calc_{E})) \notag \\
& + (\nabla^{g}_{E} (\calc_{X}) - \nabla^{g}_{X} (\calc_{E}))\calc_{Y} +
\calc_{X} \calc_{Y} \notag \\
=& -\nabla^{g}_{X \circ Y}(\calc_{E}) - \calc_{X} \nabla^{g}_{Y}(\calc_{E})
- \nabla^{g}_{X}(\calc_{E}) \calc_{Y} + \calc_{X} \calc_{Y} \notag \\
=& \;[\calc_{X \circ Y}, \calv] - \calc_{X \circ Y} +
\calc_{X}[\calc_{Y}, \calv]- \calc_{X} \calc_{Y} + [\calc_{X},
\calv] \calc_{Y} - \hspace{-1pt}\calc_{X} \calc_{Y} + \calc_{X} \calc_{Y} \notag
\\
=&\; 0 .\notag
\end{align}
\end{proof}

\newpage
\subsection{The general case} $ $\\
\begin{theorem}\label{U:General}
Let $((M,0),D,K, \dr,\calc, \calu, \calv, g)$ be a Frobenius type
structure with logarithmic pole along $D$, $\xi \in \calo(K)$
and $\dr (\xi_{\mid M \setminus D}) = 0$ with
\begin{align}
& (IC) \quad \calc_{\bullet} \xi\hspace{-2pt}\mid_0 : \dmld_{0} \longrightarrow K_{0}
\quad \text{injective,} \notag \\
& (GC) \quad \xi\hspace{-2pt}\mid_0 \; \text{and its images under iteration of the maps
} \; \calc_{X}, X \in \dmld_{0} \;  \notag  \\
& \qquad \quad \; \text{and} \;\; \calu: K_{0} \longrightarrow K_{0} \;
\text{generate}\ \; K_{0},
\notag \\
& (EC) \quad \calv \xi = \frac{d}{2} \xi \quad \text{for} \; d \in
\mc \; .\notag
\end{align}
Then up to canonical isomorphism there exist unique data $((\tilde{M},0),
\tilde{D}, \circ, e ,E, \tilde{g},$ $i,j)$ with $((\tilde{M},0),
\tilde{D}, \circ, e ,E, \tilde{g})$ a Frobenius manifold with
logarithmic pole along $\tilde{D}$,\\ $i: (M,0) \rightarrow
(\tilde{M},0)$ is an embedding such that $i(M) \cap \tilde{D} = i(D)$ and $j: K \rightarrow
Der_{\tilde{M}}(log \tilde{D})_{\mid i(M)}$ is an isomorphism above
$i$ of germs of vector bundles which identifies the logarithmic
Frobenius type structure on $K$ with the natural logarithmic Frobenius type
structure on $Der_{\tilde{M}}(log \tilde{D})\hspace{-2pt}\mid_{i(M)}$.
\end{theorem}$ $\\
Here we also have to demand the existence of a \textbf{section} $\xi$ which is flat on $M\setminus D$
and fulfills $(IC), (GC)$ and $(EC)$ unlike in \cite{HM} where they could extend a vector $\xi$ to a flat
section on $M$.\\[7pt]

\noindent First we give an overview of the proof.\\[7pt]
\textbf{(A)} Proposition $\ref{FTSlogDoneone}$ is used to get a
$(logD-trTLEP(w))$-structure. Then we prove lemma $\ref{1:longlem}$
below. This lemma together with lemma \ref{1:Punfold} gives us an unfolding of the
$(logD-trTLEP(w))$-structure which fulfills property $(\ref{1:unfoldeq})$ below for
given $f_1 , \ldots , f_n
\in \calo_{M \times \mc^{l}}$.\\[7pt]
\textbf{(B)} We investigate transformation properties of elementary
sections with respect to pullback in the context of lemma
$\ref{1:longlem}$. This is needed in order to prove uniqueness of the constructed Frobenius
manifold.\\[7pt]
\textbf{(C)} We choose an unfolding of the given
$(logD-trTLEP(w))$-structure so that the map
\[
\tilde{\calc}_{\bullet} \xi \mid_0 : Der_{\tilde{M}}(log\tilde{D}) \longrightarrow K_0
\]
is an isomorphism which is possible by lemma $\ref{1:longlem}$. We
prove that an unfolding with this property is a universal unfolding.
But because two universal unfoldings are isomorphic, this shows that
the property above characterizes the unfolding up to isomorphism.
Now an application of proposition $\ref{1:isocase}$ shows that this
unfolding gives rise to a Frobenius manifold which is therefore also unique up to
isomorphism.\\[5pt]

\noindent (\textbf{A})\\
Let $((M,0),D, H, \nabla, P)$ be a $(logD - trTLEP(w))$-structure. We want to choose \textbf{global}
sections $(v_1,\ldots ,v_n)$ which are a basis of $H$ so that the connection matrix has a convenient form for proving
the existence of an unfolding. Let $\nabla^{res}$ be the residue connection on $H_{\mid\{ \infty \} \times (M \setminus D,0)}$.
We choose the sections so that their restrictions to $\{ \infty \} \times (M,0)$ are elementary sections with respect to the
residue connection. We call such a choice of sections adapted to $((M,0),D, H, \nabla, P)$.
Now let $((M \times \mc^l,0),D\times \mc^l, \tilde{H}, \tilde{\nabla})$ be an unfolding of the above structure. Then there exist unique adapted sections $(\tilde{v}_1,\ldots ,\tilde{v}_n)$ of $((M \times \mc^l,0),D\times \mc^l, \tilde{H}, \tilde{\nabla})$ so that $(\tilde{v}_1\hspace{-3pt}\mid_{\mathbb{P}^{1} \times M \times \{ 0 \}}, \ldots , \tilde{v}_n\hspace{-3pt}\mid_{\mathbb{P}^{1} \times M \times \{ 0\}}) = (v_1, \ldots , v_n)$. We call these sections a canonical extension of $(v_1, \ldots ,v_n)$.\\
Recall the proof of proposition \ref{FTSlogDoneone}. There we observed that $K = H_{ \mid \{0\}  \times (M,0) }$ is canonically
isomorphic to $H_{\mid \{ \infty \} \times (M,0)}$. Therefore it is possible to shift structure from one bundle to the other. This
is the case for the residue connection $\nabla^{res}$ which when shifted to $K$ can be identified with $\nabla^r$. Consider a global section $v$ of $H$. Then its restriction to $\{ 0 \} \times (M,0)$ is $\nabla^r$-flat if and only if its restriction to $\{\infty \} \times (M,0)$ is $\nabla^{res}$-flat.
\begin{lemma}\label{1:longlem}
Let $w \in \mathbb{Z}$ and $((M,0),D, H, \nabla, P)$ be the $(logD -
trTLEP(w))$-structure which fulfills the conditions above. Choose
adapted sections $(v_{1}, \ldots, v_{n})$ of H such that $v_{1 \mid \{ 0 \} \times (M,0)}
= \xi$.
Choose $l \in \mathbb{N}$ and $n$ functions $f_{1}, \ldots, f_{n}
\in \mathcal{O}_{(M \times \mc^{l},0)}$ with $f_{i}\mid_{(M \times
\{0\},0)} = 0$. Let $(t_{1}, \ldots t_{m},y_{1}, \ldots ,y_{l}) =
(t,y)$ be coordinates on $(M \times \mc^{l},0)$.\\
Then there exists a unique unfolding $(\tilde{H} \rightarrow
\mathbb{P}^{1} \times (M \times \mc^{l}, 0), \tilde{\nabla})$ of
$(H, \nabla)$ with the following properties. $\tilde{H}$ is a trivial vector bundle with
flat connection over $(\mc^{*} \times ((M \setminus D) \times \mc^{l},
0))$ with logarithmic poles along $\{\infty\} \times ( M \times
\mc^{l},0 ) \cup \mathbb{P}^1\setminus \{ 0 \} \times (D \times \mc^l ,0)$ and a pole of Poincar\'{e} rank 1 along
$\{0\} \times (M \times \mc^{l},0)$. Its restriction to
$\mathbb{P}^{1} \times (M \times \{0\},0)$ is $(H, \nabla)$. Let
$(\tilde{v}_{1}, \ldots , \tilde{v}_{n})$ be the canonical extensions
of $(v_{1}, \ldots , v_{n})$. Define $\tilde{K}:= \tilde{H}_{\mid
\{0\} \times (M \times \mc^{l}, 0)}$ and $ \tilde{\calc},
\tilde{\calu}$ as above. Then
\begin{equation}\label{1:unfoldeq}
\tilde{\calc}_{\partial_{y_{\alpha}}} \tilde{v}_{1} = \sum^{n}_{i=1}
\frac{\partial f_{i}}{\partial y_{\alpha}} \tilde{v}_{i} \quad
\text{for} \quad \alpha = 1, \ldots , l \; .
\end{equation}
\end{lemma}$ $\\
\begin{proof}
Suppose for a moment that $(\tilde{H}, \tilde{\nabla})$ were already
constructed. The connection matrix $ \Omega$ with respect to the
basis $(\tilde{v}_{1}, \ldots, \tilde{v}_{n})$,
\[
\tilde{\nabla}(\tilde{v}_1, \ldots, \tilde{v}_n) = (\tilde{v}_1, \ldots , \tilde{v}_n) \cdot \Omega,
\]
would take the form
\[
\Omega = \sum_{i=1}^{r} A_{i} \frac{dt_{i}}{t_{i}} + \frac{1}{z}
\sum^{r}_{i=1} C_{i} \frac{dt_{i}}{t_{i}} + \frac{1}{z} \sum_{k =
r+1}^{m} C_{k} dt_{k} + \frac{1}{z} \sum_{\alpha =1}^{l} F_{\alpha}
dy_{\alpha} + (\frac{1}{z^{2}} U + \frac{1}{z} V )
dz
\]
with matrices
\[
A_{i},C_{i},C_{j},F_{\alpha},U, V, \in M(n \times n, \calo_{M
\times \mc^{l}})\; .
\]
This follows from $\tilde{v}_{i}$ elementary and from the pole
orders of $(\tilde{H}, \tilde{\nabla})$ along $\{0, \infty\} \times
( M \times \mc^{l},0)$.\\
\noindent In the following $i,j \in
\{1, \ldots , r \},\; k,l \in \{ r+1, \ldots m \}$. The flatness
condition $d \Omega + \Omega \wedge \Omega =0$ implies
\begin{align}
\frac{\p A_{i}}{\p t_{k}} &= 0 \label{C:0}, \\
\frac{\partial A_{i}}{\partial y_{\alpha}} &=  0, \label{C:1} \\
[A_{i},A_{j}] &= - t_{i}\frac{\partial A_{j}}{\partial t_{i}} + t_{j}
\frac{\partial A_{i}}{\partial t_{j}}, \label{C:2} \\
[C_{i},C_{j}] &=
[C_{i},C_{k}] = [C_{k},C_{l}] = 0, \label{C:3} \\
[C_{i} , F_{\alpha}] &= [C_{k},F_{\alpha}] = 0, \label{C:4} \\
[F_{\alpha},F_{\beta}] &= 0, \label{C:5} \\
t_{j} \frac{\p C_{i}}{\p t_{j}} + [A_{j}, C_{i}]  &= t_{i} \frac{\p
C_{j}}{\p t_{i}} + [A_{i}, C_{j}], \label{C:6} \\
\frac{\p C_{i}}{\p t_{k}} &= t_{i} \frac{\p C_{k}}{\p t_{i}} +
[A_{i},C_{k}], \label{C:7} \\
\frac{\p C_{l}}{\p t_{k}} &= \frac{\p C_{k}}{\p t_{l}}, \label{C:8} \\
\frac{ \p C_{i}}{\p y_{\alpha}} &= t_{i} \frac{\p F_{\alpha}}{\p
t_{i}} + [A_{i}, F_{\alpha}], \label{C:9} \\
\frac{\p C_{k}}{\p y_{\alpha}} &= \frac{\p F_{\alpha}}{\p t_{k}}, \label{C:10} \\
\frac{\p F_{\alpha}}{\p y_{\beta}} &= \frac{\p
F_{\beta}}{\p y_{\alpha}}, \label{C:11} \\
[U, C_{j}] &= [U, C_{k}] = 0, \label{C:12} \\
[U, F_{\alpha}] &= 0, \label{C:13} \\
t_i \frac{\p U}{\p t_{i}} &= [V, C_{i}] -
C_{i} + [U, A_{i}], \label{C:14} \\
\frac{\p U}{\p t_{k}} &= [V, C_{k}] - C_{k}, \label{C:15} \\
\frac{\p U}{\p y_{\alpha}} &= [V,F_{\alpha}] - F_{\alpha}, \label{C:16} \\
t_i \frac{\p V}{\p t_{i}} &= [ V, A_{i}], \label{C:20}\\
\frac{\p V}{\p t_{k}} &= 0 , \label{C:21}\\
\frac{\p V}{\p y_{\alpha}} &= 0,  \label{C:22} \\
& \notag \\
& \hspace{-156pt}\text{The condition (\ref{1:unfoldeq}) would mean} \notag \\
& \notag \\
(F_{\alpha})_{i1} &= \frac{\p f_{i}}{\p y_{\alpha}}. \label{C:23}
\end{align}
The proof consists of three parts. In parts (I) and (II) we restrict
to the case $l = \alpha = 1$. In part (I) we show inductively
uniqueness and existence of matrices $C_{i},C_{j},F_{1}, U, V,W$
with (\ref{C:3})-(\ref{C:23}) and with coefficients in
$\calo_{M,0}[[y_{1}]]$. In part (II) their convergence will be
proved with the Cauchy-Kovalevski theorem. In part (III) the general
case will be proved by induction in $l$.
\paragraph{Part (I)} Suppose
$l = \alpha = 1$ and $y_{1} = y$. Define for $ w \in
\mathbb{Z}_{\geq 0}$
\begin{align}
\calo_{M,0}[y]_{\leq w} &:= \sum_{k=0}^{w} \calo_{M,0} \cdot y^{k},
\notag \\
\calo_{M,0}[y]_{> w} &:= \calo_{M,0}[y] \cdot y^{w+1}, \notag \\
\calo_{M,0}[[y]]_{> w} &:= \calo_{M,0}[[y]] \cdot y^{w+1} \notag
\end{align}
and
\begin{align}
M(w) &:= M(n \times n,\calo_{M,0} \cdot y^{k}), \notag \\
M(> w) &:= M(n \times n,\calo_{M,0}[y]_{> w}), \label{1:Mgw} \\
M( \leq w) &:= M(n \times n, \calo_{M,0}[y]_{\leq w}). \notag
\end{align}
$ $\\[10pt]
\textit{Beginning of the induction for} $w=0$: The connection
matrix $\Omega^{(0)}$ of $(H,\nabla)$ with respect to the basis
$v_{1}, \ldots, v_{n}$ takes the form
\[
\Omega^{(0)} = \sum_{i=1}^{r} A_{i} \frac{dt_{i}}{t_{i}} +
\frac{1}{z} \sum^{r}_{i=1} C_{i}^{(0)} \frac{dt_{i}}{t_{i}} +
\frac{1}{z} \sum_{j = r+1}^{m} C_{j}^{(0)} dt_{j} + (\frac{1}{z^{2}}
U^{(0)} + \frac{1}{z} V^{(0)} ) dz
\]
with matrices $A_{i},
C_{i}^{(0)},C_{j}^{(0)},U^{(0)},V^{(0)} \in M(0)$. The
flatness condition $d\Omega^{(0)} + \Omega^{(0)} \wedge \Omega^{(0)}
= 0$
is equivalent to the equations
$(\ref{C:0}),(\ref{C:2}),(\ref{C:3}),(\ref{C:6})$, $(\ref{C:7})$, $(\ref{C:8})$, $(\ref{C:12})$,
$(\ref{C:14}),(\ref{C:15}),$
$(\ref{C:20}),(\ref{C:21})$ for $w=0$.
The equations
$(\ref{C:1}),(\ref{C:4})$, $(\ref{C:5})$, $(\ref{C:9}),(\ref{C:10}),(\ref{C:11}),(\ref{C:13}),$ $(\ref{C:16}),(\ref{C:22})$
are empty.\\[15pt]
\textit{Induction hypothesis for} $w \in \mathbb{Z}_{\geq 0}$:
Unique matrices $C_{i}^{(k)},C_{j}^{(k)},U^{(k)},V^{(k)} \in
M(k)$ for $0 \leq k \leq w$ and $F_{1}^{(k)} \in M(k)$ for $0 \leq k
\leq w-1$ are constructed such that the matrices
\[
C_{i}^{(\leq w)} := \sum_{k=0}^{w} C_{i}^{(k)} \in M(\leq w)\; ,
\]
the analogously defined matrices $C_{j}^{(\leq w)}, U^{(\leq
w)},V^{(\leq w)}, \hspace{-1pt}F_{1}^{(\leq w-1)}$ and the matrices $A_i$ satisfy
$(\ref{C:3}),(\ref{C:6}), (\ref{C:7}), (\ref{C:8}),(\ref{C:12})
,(\ref{C:14}),
(\ref{C:15}),(\ref{C:20}),(\ref{C:21})$ modulo $M(>w)$,
$(\ref{C:4}),(\ref{C:5}),$ $(\ref{C:9}),(\ref{C:10}),(\ref{C:13}),(\ref{C:16}),(\ref{C:22})$
modulo $M(>w-1)$ and $(\ref{C:23})$ modulo
$\calo_{M,0}[[y]]_{>w-1}$.\\
\textit{Induction step from $w$ to $w+1$}: It consists of three
steps:\\
\begin{enumerate}
\item Construction of a matrix $F_{1}^{(w)} \in M(w)$ such that the
matrix $F_{1}^{(\leq w)} = F_{1}^{(\leq w-1)} + F_{1}^{(w)}$
together with the matrices $C_{i}^{(\leq w)},C_{j}^{(\leq
w)},U^{(\leq w)}, V^{(\leq w)}$ satisfies
$(\ref{C:4}),(\ref{C:13})$ mod $M(>w)$ and $(\ref{C:23})$ mod $\calo_{M,0}[[y]]_{> w} .$ \\
\item Construction of matrices
$C_{i}^{(w+1)}, C_{k}^{(w+1)}, U^{(w+1)}, V^{(w+1)} \in
M(w+1)$ such that the matrices $C_{i}^{(\leq w+1)} \hspace{-1pt}= C_{i}^{(\leq
w)} + C_{i}^{(w+1)}$ and the analogously defined matrices
$C_{j}^{(\leq w+1)},U^{(\leq w+1)},$ $V^{(\leq w+1)}$
and the matrix $F_{1}^{(\leq w)}$ satisfy
$(\ref{C:9}),(\ref{C:10})$, $(\ref{C:16}),(\ref{C:22})$
mod $M(>w)$.\\
\item Proof of
$(\ref{C:3}),(\ref{C:6}),(\ref{C:7}),(\ref{C:8}),(\ref{C:12}),(\ref{C:14}),(\ref{C:15}),(\ref{C:20})$,
$(\ref{C:21})$ mod $M(>w+1)$.
\end{enumerate}$ $\\
\textbf{(1)} The matrices $C_{i}^{(\leq w)},C_{k}^{(\leq w)}$ and
$U^{(\leq w)}$ generate an algebra of commuting matrices in $M(n
\times n, \calo_{M,0}[y]/\calo_{M,0}[y]_{>w})$. Because of  the
generation condition (GC), the image of the column vector
$(1,0,\ldots,0)^{tr}$ under the action of this algebra is the whole
space $M(n \times 1,\calo_{M,0}[y]/\calo_{M,0}[y]_{>w})$. This shows
two things:
\begin{itemize}
\item This algebra contains matrices $E_{i}^{(\leq w)} \in M(n \times n, \calo_{M,0}[y]_{\leq w})$ for any $i = 1, \ldots, n$ with the first column
\[
(E_{i}^{(\leq w)})_{j1} = \delta_{ij} .
\]
\item Any matrix in $M(n \times n,\calo_{M,0}[y]_{\leq w})$ which
commutes with the matrices $C_{i}^{(\leq w)},C_{j}^{(\leq w)}$ and
$U^{(\leq w)}$ mod $M(>w)$ is modulo $M(>w)$ a linear combination of
the matrices $E_{i}^{(\leq w)}$ with coefficients in
$\calo_{M,0}[y]_{\leq w}$.
\end{itemize}
Therefore the matrix $F_{1}^{(\leq w)} \in M(\leq w)$ which is
defined by
\begin{equation}
F_{1}^{(\leq w)} = \sum_{i=1}^{n} \frac{\p f_{i}}{\p y} \cdot
E_{i}^{(\leq w)} \quad \text{mod}\;\;  M(n \times n,
\calo_{M,0}[[y]]_{>w}) \label{C:24}
\end{equation}
is the unique matrix which satisfies $(\ref{C:4}),(\ref{C:13})$ mod
$M(>w)$, $(\ref{C:23})$ mod\\ $\calo_{M,0}[[y]]_{>w}$ and
$F_{1}^{(\leq w)} = F_{1}^{(\leq
w-1)} + F_{1}^{(w)}$ for some $F_{1}^{(w)} \in M(w)$.\\[14pt]
\textbf{(2)}
\begin{align}
& \frac{\p C_{i}}{\p y} = t_{i} \cdot \frac{\p F_{1}}{\p
t_{i}} + [A_{j},F_{1}] \quad mod \quad M(>w) \notag \\
& \Rightarrow C_{i}^{(\leq w+1)} \quad \text{constructed}, \notag
\end{align}
and so for the other matrices $C_j^{(\leq w+1)},U^{(\leq w+1)}, V^{(\leq w+1)}$ with the corresponding equations.\\
\noindent \textbf{(3)} Now we check that the derivatives by $\frac{\p}{\p y}$
of the remaining equations hold modulo $M(>w+1)$. Because of the derivative
we calculate modulo $M(>w)$ and therefore we can use all equations modulo $M(>w)$. For example one calculates:\\
(\ref{C:3})
\begin{align}
\frac{\p}{\p y} \left[C_{i},C_{j}\right]
\hspace{2 pt} &\hspace{-3 pt} \overset{(\ref{C:9})}{=}
\left[t_i \frac{\p F_{1}}{\p t_{i}}+ \left[A_{i},F_{1}
\right], C_{j}\right] +
\left[C_{i},t_j \frac{\p F_{1}}{\p t_{j}} +
\left[A_{j},F_{1}\right]\right] \notag \\
&\hspace{4 pt}= \left[t_i \frac{\p F_{1}}{\p t_{i}},
C_{j}\right] + \left[\left[A_{i},F_{1}
\right],C_{j}\right] + \left[C_{i},
t_j \frac{\p F_{1}}{\p t_{j}}\right] + \left[C_{i},
\left[A_{j},F_{1} \right]\right] \notag \\
&\hspace{-2 pt}\overset{(\ref{C:4})}{=} -\left[F_{1}
,t_{i} \frac{\p C_{j}}{\p t_{i}}\right] + \left[F_{1} ,\left[C_{j},A_{i}\right]\right] +\left[F_{1}
,t_{j}\frac{\p C_{i}}{\p t_{j}}\right] -  \left[F_{1} ,\left[C_{i},A_{j}\right]\right] \notag
\\
&\hspace{-2 pt}\overset{(\ref{C:6})}{=} 0\; . \notag
\end{align}
The other equations are similar or easier. This finishes the proof of the induction step from $w$ to $w+1$. It shows uniqueness and existence of matrices $C_i, C_j,F_1 ,U$, $V,W \in M(n \times n, \calo_{M,0}[[y]])$ with $(\ref{C:3}) - (\ref{C:23})$ and with restrictions\\ $(C_i,C_j,U,V)\mid_{y=0} = (C_i^{(0)},C_j^{(0)},U^{(0)},V^{(0)})$.\\
\noindent \paragraph{Part (II)} Now we have to show holomorphy of these matrices. We
want to apply the Cauchy-Kovalevski theorem in the following
form(\cite{Fo}(1.31),(1.40),(1.41); there the setting is real
analytic, but proofs and statements hold also in the complex
analytic setting):\\ Given $N\in \mathbb{N}$ and matrices $H_i, L \in
M(N \times N, \mc\{t_1, \ldots , t_m,y,x_1, \ldots , x_N\})$ there
exists a unique vector $\Phi \in M(N \times 1, \mc\{t_1, \ldots ,
t_m,y\})$ with
\begin{align}
&\frac{\p \Phi}{\p y} = \sum_{i=1}^{m} H_i(t,y,\Phi) \frac{\p \Phi}{\p t_i} + L(t,y,\Phi), \label{CKS:1} \\
&\Phi(t,0) = 0. \label{CKS:2}
\end{align}
We will construct a system (\ref{CKS:1}) - (\ref{CKS:2}) with $N
=(m+3)n^2$ such that it will be satisfied with the entries of the
matrices $C_i - C_i^{(0)}, U- U^{(0)},V - V^{(0)}$ as
entries of $\Phi$. The system will be built from the following
equations
\begin{align}
\frac{\p C_{i} -C_{i}^{(0)}}{\p y} &= t_{i} \frac{\p F_{1}}{\p
t_{i}}
+[A_{i},F_{1}], \label{CKeq3} \\
\frac{\p C_{k} - C_{k}^{(0)}}{\p y} &= \frac{\p F_{1}}{\p t_{k}},
\label{CKeq4} \\
\frac{\p U - U^{(0)}}{\p y} &= [V,F_{1}] -F_{1}, \label{CKeq5} \\
\frac{\p V - V^{(0)}}{\p y} &= 0 \label{CKeq7}
\end{align}
and equations $(\ref{CKeq1}),(\ref{CKeq2})$ with which one can
express the entries of $F_1$ as functions of the entries of $\Phi$.
$ $\\
The commutative subalgebra of $M(n \times n, \calo_{M,0}[[y]])$,
which is generated by the matrices $C_{1}, \ldots , C_{m},U$, is a
free $\calo_{M,0}[[y]]$-module of rank n. Choose monomials $G^{(j)},
j = 1, \ldots ,n$, in the matrices $C_{1},\ldots, C_{m},U$ which
form an $\calo_{M,0}[[y]]$-basis of this module. Then the matrix
$(G_{i1}^{(j)})_{ij}$ of the first columns of the matrices $G^{(j)}$
is invertible in $M(n \times n, \calo_{M,0}[[y]])$. Equation
(\ref{C:24}) gives
\begin{equation}\label{CKeq1}
 F_{1} = \sum_{j=1}^{n} g_{j} G^{(j)}
\end{equation}
with coefficients $g_{j} \in \calo_{M,0}[[y]]$ such that
\begin{equation}\label{CKeq2}
( \frac{\p f_{1}}{\p y}, \ldots , \frac{\p f_{n}}{\p y})^{tr} =
(G_{i1}^{(j)}) \cdot (g_{1}, \ldots , g_{n})^{tr}.
\end{equation}
Replacing the entries of the matrices $C_{i} - C_{i}^{(0)},C_{k} -
C_{k}^{(0)}, U - U^{(0)},V - V^{(0)}$ by indeterminates
$x_{1}, \ldots x_{N}$, the coefficients of the matrices $G^{(j)}$
become elements of $\mc\{t\}[x_{1},\ldots,x_{N}]$, and the
$g_{j}$ become elements of
$\mc\{t_{1},\ldots,t_{m},y,x_{1},\ldots,x_{N}\}$. One obtains from
$(\ref{CKeq3}) - (\ref{CKeq2})$ a system
$(\ref{CKS:1}),(\ref{CKS:2})$. Now the theorem of Cauchy-Kovalevski
shows $C_{i},C_{k},F_{1},U,V \in M(n \times n,
\calo_{M \times \mc, 0})$. This shows lemma \ref{1:longlem} in the case of $l=1$.\\[-3pt]
\paragraph{Part (III)} By induction in $l$
one obtains a slightly weaker version of the lemma, namely with
formula (\ref{1:unfoldeq}) replaced by
\[
(\tilde{\calc}_{\p y_{\alpha}} \tilde{v}_{1})\mid_{{y_{\alpha + 1}=
\ldots =y_{l}= 0}} = (\sum_{i=1}^{n} \frac{\p f_{i}}{\p y_{\alpha}}
\tilde{v}_{i}) \mid_{{y_{\alpha + 1}= \ldots =y_{l}= 0}}
\]
for $\alpha = 1, \ldots ,l$. This equation is equivalent to
$(\ref{C:23})$ with the same restrictions. But now one has a
connection matrix as in the third line of the proof of lemma
\ref{1:longlem}. Flatness gives $(\ref{C:0}) - (\ref{C:22})$ . The
equation $(\ref{C:11})$  together with the formula above gives the
formula at the end of lemma $\ref{1:longlem}$.
\end{proof}
\vspace{5pt}
\noindent It rests to prove that the pairing $P$ extends to the unfolding. We need a lemma about the rigidity of logarithmic poles.
\begin{lemma}
Let $(L' \rightarrow (\mc^{\ast},0) \times (M,0), \nabla)$ be the germ of a holomorphic vector bundle with flat connection $\nabla$, and let $(L^{(0)} \rightarrow (\mc,0) \times \{0 \}, \nabla)$ be an extension of $(L',\nabla)\mid_{(\mc^{\ast},0) \times \{ 0 \}}$ with a logarithmic pole at $0$. Then an extension $(L \rightarrow (\mc,0) \times (M,0), \nabla)$ of $(L', \nabla)$ with a logarithmic pole along $\{ 0 \} \times (M,0)$ exists with $(L, \nabla)\mid_{(\mc,0) \times \{0 \}} = (L^{(0)}, \nabla)$. It is unique up to canonical isomorphism and it is isomorphic to the pullback $\varphi^{\ast}(L^{(0)}, \nabla)$ where $\varphi: (M,0) \rightarrow \{0 \}$.
\end{lemma}
\begin{proof}
See \cite{Sab} III.1.20 .
\end{proof}
\begin{lemma}\label{1:Punfold}
Let $((M \times \mc^l,0), D \times \mc^l  , \tilde{H},\tilde{\nabla})$ be the unfolding of the $(logD-trTLEP(w))$-structure $(H,\nabla, P)$ given above. Then $P$ extends to $\calo(\tilde{H})$ and $((M \times \mc^l,0), \tilde{D}, \tilde{H},\tilde{\nabla},P)$ is a $log\tilde{D}-trTLEP(w)$-structure.
\end{lemma}
\begin{proof}
It is sufficient to consider an unfolding in one parameter $y$. For some representative of $\tilde{H}$ the pairing $P$ extends to a $\nabla$-flat pairing on the restriction to $(\mc^{\ast},0) \times (M\setminus D \times \mc,0)$. We have to show that it takes values on $\calo(\tilde{H})$ in $z^w \calo_{\mathbb{P}^1 \times M \times \mc}$. A priori the values are in $\calo_{\mc^{\ast} \times M \setminus D \times \mc}$. Recall that $\nabla$ has a logarithmic pole along $(\mc^{\ast},0) \times (D,0)$. Because of the rigidity of logarithmic poles the extension to $(\tilde{H},\tilde{\nabla})\mid_{\mc^{\ast}\times M\setminus D_{sing} \times \mc}$ is trivial, therefore $P$ takes values in $\calo_{\mc^{\ast}\times M\setminus D_{sing} \times \mc}$. A codimension $2$ argument shows that $P$ takes values $\calo_{\mc^{\ast}\times M  \times \mc}$. The same reasoning shows that on $\mathbb{P}^{1}\setminus \{0\} \times M \times \mc$ the pairing $P$ takes values in $z^w\calo_{\mathbb{P}^{1}\setminus \{0\} \times M \times \mc}$.
It rests to prove the corresponding property on $\mc \times M \times \mc$. Denote $n:= rk H$ and let $(z,t_1, \ldots, t_m,y)$ be coordinates on $(\mc \times M \times \mc,0)$. Choose an $\calo_{\mc \times M \times \mc,0}$-basis $(\tilde{v}_1, \ldots , \tilde{v}_n)$ of $\calo(\tilde{H})$ with connection matrix
\[
\Omega = \frac{1}{z}\sum_{i=1}^r \check{C}_i \frac{dt_i}{t_i} + \frac{1}{z} \sum_{i=r+1}^{m} \check{C}_k dt_k + \frac{1}{z} \check{F}dy+ \frac{1}{z} \check{U} dz
\]
with matrices $\check{C}_i,\check{C}_k,\check{F},\check{U} \in M(n \times n, \calo_{\mc \times M \times \mc,0})$ and the matrix
\[
R:= ( P(\tilde{v}_i,\tilde{v}_j)) \in M(n \times n, \calo_{\mc^{\ast} \times M \times \mc,0}).
\]
Flatness and $z$-sesquilinearity of the pairing give
\[
dR(z,t,y) = \Omega^{tr}(z,t,y)R(z,t,y) + R(z,t,y)\Omega(-z,t,y),
\]
that means,
\begin{align}
z \frac{\p}{\p z} R(z,t,y) &= \frac{1}{z} \check{U}^{tr}(z,t,y)R(z,t,y) - \frac{1}{z}R(z,t,y)\check{U}(-z,t,y), \label{1:Punfold1}\\
t_i \frac{\p}{\p t_i} R(z,t,y) &= \frac{1}{z} \check{C}_i^{tr}(z,t,y)R(z,t,y) - \frac{1}{z}R(z,t,y)\check{C}_i(-z,t,y), \label{1:Punfold2}\\
\frac{\p}{\p t_k} R(z,t,y) &= \frac{1}{z}\check{C}_k^{tr}(z,t,y)R(z,t,y) - \frac{1}{z}R(z,t,y)\check{C}_k(-z,t,y), \label{1:Punfold3}\\
\frac{\p}{\p y}R(z,t,y) &= \frac{1}{z}\check{F}^{tr}(z,t,y)R(z,t,y) - \frac{1}{z}R(z,t,y)\check{F}(-z,t,y). \label{1:Punfold4}
\end{align}
Write $R$ as a power series
\[
R(z,t,y) = \sum_{l=0}^{\infty} R^{(l)}(z,t) \;\;\text{with}\;\; R^{(l)} \in M(n \times n, \calo_{\mc^{\ast}\ \times M,0} \cdot y^l)
\]
and define
\[
R^{(\leq l)}(z,t,y) := \sum_{j=0}^{l} R^{(j)}(z,t),
\]
analogously for $\check{C}_i,\check{C}_k,\check{F},\check{U}$, with $\check{C}_i^{(l)},\check{C}_k^{(l)},\check{F}^{(l)},\check{U}^{(l)} \in M(n \times n, \calo_{\mc \times M,0} \cdot y^l)$. Then
$R^{(0)} \in M(n \times n, z^w \calo_{\mc \times M,0})$ because $(H, \nabla, P)$ is a $(logD-trTLEP(w))$-structure.\\
\textit{Induction hypothesis for $k \in \mathbb{Z}_{\geq 0}$}:
\[
R^{(\leq l)} \in M(n \times n, z^w \calo_{\mc \times M \times \mc,0}).
\]
\textit{Induction step from $l$ to $l+1$}: Recall the definition of $M(>l)$ in (\ref{1:Mgw}). The equations (\ref{1:Punfold1}),(\ref{1:Punfold2}) and (\ref{1:Punfold3}) show that one has modulo $M(>l)$
\begin{align}
& \check{C}_i^{(\leq l) tr}(0,t,y)[z^{-w}R^{(\leq l)}(z,t,y)]\mid_{z=0} \notag \\
\equiv & [z^{-w}R^{(\leq l)}(z,t,y)]\mid_{z=0}\check{C}_i^{(\leq l)}(0,t,y), \notag \\
& \check{C}_k^{(\leq l) tr}(0,t,y)[z^{-w}R^{(\leq l)}(z,t,y)]\mid_{z=0} \notag \\
\equiv & [z^{-w}R^{(\leq l)}(z,t,y)]\mid_{z=0}\check{C}_k^{(\leq l)}(0,t,y), \notag \\
& \check{U}^{(\leq l) tr}(0,t,y)[z^{-w}R^{(\leq l)}(z,t,y)]\mid_{z=0} \notag \\
\equiv & [z^{-w}R^{(\leq l)}(z,t,y)]\mid_{z=0}\check{U}^{(\leq l)}(0,t,y), \notag
\end{align}
Because of the generation condition (GC) the matrix $\check{F}^{(\leq l)}(0,t,y)$ is an element of the commutative subalgebra of
$M(n\times n, \calo_{M,0}[y])/M(>l)$ which is generated by $\check{C}_1^{(\leq l)}, \ldots , \check{C}_m^{(\leq l)}, \check{U}^{(\leq l)}$. Therefore modulo
$M(>l)$
\begin{align}
&\check{F}^{(\leq l)tr}(0,t,y)[z^{-w}R^{(\leq l)}(z,t,y)]\mid_{z=0} \notag \\
\equiv &[z^{-w}R^{(\leq l)}(z,t,y)]\mid_{z=0}\check{F}^{(\leq l)}(0,t,y). \notag
\end{align}
This together with (\ref{1:Punfold4}) completes the induction step.
\end{proof}
\noindent \textbf{(B)} Transformation properties \\

In this section we want to clarify the transformation properties of the global sections
$v_1, \ldots , v_n$, which were defined in lemma $\ref{1:longlem}$, under pullback by a map
$\varphi : M \times \mc^{l'} \rightarrow M \times \mc^{l}$ of a given $log(D)-trTLEP(w)$-structure.
Furthermore the meaning of the base change matrix $B$ between
the elementary sections $s_{j}$ around the divisor $\lbrace \frac{1}{z} \cdot
t_{1} \cdot \ldots \cdot t_{r} = 0 \rbrace$ and the global sections
$v_{j}$ of an adapted basis will be investigated.\\
W.l.o.g. we can assume that
\[
 B = I + B^{(1)}z^{-1} + B^{(2)} z^{-2} + \ldots \;.
\]
Let $rank\; H =n$ and $s_{1},\ldots,s_{n}$ be elementary sections as
above.\\ There exists
$\alpha_{j},\beta_{j}^{(1)},\ldots,\beta_{j}^{(r)}\in \mc$ with $j
\in \lbrace 1, \ldots , n \rbrace$,
\begin{align}
A_{j} &\in H^{\infty}(e^{-2 \pi i \alpha_{j}}, e^{- 2\pi i
\beta_{j}^{(1)}}, \ldots, e^{- 2 \pi i \beta_{j}^{(r)}}), \notag \\
s_{j} &= z^{\alpha_{j} - \frac{N^{(0)}}{2 \pi i}} \prod_{k=1}^{r}
t_{k}^{\beta_{j}^{(k)} - \frac{N^{(k)}}{2 \pi i}} A_{j} \notag
\end{align}
and for $k \in \{ 1, \ldots ,r \}$
\begin{align}
M^{(0)} &:= \text{Monodromy around}\;\; z=0 ,\quad &&N^{(0)} =
\log(\text{unipotent part of}\;\; M^{(0)}), \notag \\
M^{(k)} &:= \text{Monodromy around}\;\; t_{k}=0 ,\quad &&N^{(k)} =
\log(\text{unipotent part of}\;\; M^{(k)}). \notag
\end{align}
With respect to the basis $A_1, \ldots,  A_n$ we have matrices $M^{(k)}_{mat},
M^{(k)}_{s,mat}, M^{(k)}_{u,mat},$ $N^{(k)}_{mat}$ with
\[
M^{(k)}_{s,mat} = \left( \begin{matrix} e^{-2 \pi i \beta_{1}^{(k)}} & & \\
& \ddots & \\ & & e^{-2 \pi i \beta^{(k)}_{n}}
\end{matrix}\right).
\]
for $k \in \{0,1, \ldots ,r \}$ with $\beta^{(0)}_j := \alpha_j$. Observe that these matrices commute. But in general the matrix
\[
\beta_{mat}^{(k)} = \left( \begin{matrix} \beta_{1}^{(k)} & & \\
& \ddots & \\ & & \beta^{(k)}_{n}
\end{matrix}\right)
\]
does not commute with $M^{(k)}_{u,mat}$ and $N^{(k)}_{mat}$. Set
\[
\Delta := \left(\begin{matrix} z^{\alpha_{1}}\prod_{k=1}^{r}
t_{k}^{\beta_{1}^{(k)}} & & \\ & \ddots & \\ & &
z^{\alpha_{n}}\prod_{k=1}^{r} t_{k}^{\beta_{n}^{(k)}}\end{matrix}
\right).
\]
We have\vspace{10pt}
\begin{align}
\nabla_{t_{k} \ptk} s_{j} &= 0 \quad &&\text{for} \; k \in \lbrace r+1, \ldots , m \rbrace, \notag \\
\nabla_{t_{k} \ptk} s_{j} &= t_{k} \ptk (z^{\alpha_{j} -
\frac{N^{(0)}}{2 \pi i}}\prod_{l=1}^{r}t_{l}^{\beta_{j}^{(l)} -
\frac{N^{(l)}}{2 \pi i}} A_{j}) \quad &&\text{for}\; k \in \lbrace
1,
\ldots , r \rbrace \notag \\
&= \beta_{j}^{(k)} \cdot s_{j} + \frac{-1}{2 \pi i} \sum_{i =1}^{n}
\Delta_{ii}^{-1}\, (N^{(k)}_{mat})_{ij} \Delta_{jj} s_{i} . \notag
\end{align}
This means in matrix notation (for $k \in \lbrace 1, \ldots, r
\rbrace)$
\begin{equation}\label{1:Nmat}
t_{k} \ptk (s_{1}, \ldots, s_{n}) = (s_{1}, \ldots , s_{n}) \cdot \left( \beta^{(k)}_{mat} +
\frac{-1}{2 \pi i }\, \Delta^{-1} N^{(k)}_{mat}\, \Delta \right) .
\end{equation}
We have
\begin{align}
(\Delta^{-1}\, N^{(k)}_{mat} \Delta)_{ij} &=
\Delta_{ii}^{-1} (N^{(k)}_{mat})_{ij}\, \Delta_{jj} \notag \\
&= (N^{(k)}_{mat})_{ij} \cdot z^{\alpha_{j} - \alpha_{i}} \prod
_{l=1}^{r} t_{l}^{\beta_{j}^{(l)} - \beta_{i}^{(l)}}. \notag
\end{align}
Because $(s_{1}, \ldots , s_{n})$ is a holomorphic basis of
$H\hspace{-3pt}\mid_{(\mathbb{P}^{1} \setminus \lbrace 0 \rbrace)
\times M}$  we have the following restrictions on the entries of
$N^{(k)}_{mat}$:
\[
(N^{(k)}_{mat})_{ij} \neq 0\; \Longrightarrow\; \alpha_{j} \leq
\alpha_{i} \; \text{and} \; \beta_{j}^{(l)} \geq \beta_{i}^{(l)}.
\]
Now we decompose $N^{(k)}_{mat}$ in terms of powers of $z^{-1}$.
Write
\begin{align}
N^{(k)}_{mat} &= \sum_{l \geq 0} N^{(k,l)}_{mat}, \notag \\
(N^{(k,l)}_{mat})_{ij} &= \begin{cases}0 & \alpha_{j} - \alpha_{i}
\neq -l\\
(N^{(k)}_{mat})_{ij} & \alpha_{j} - \alpha_{i} = -l
\end{cases}. \notag
\end{align}
As above we have a basis $v_{1}, \ldots , v_{n}$ of global sections
such that
\[
(v_{1}, \ldots , v_{n}) = (s_{1}, \ldots, s_{n}) \cdot (\sum_{k=0}^{\infty} z^{-k}
B^{(k)}(t)) .
\]
Set
\[
\alpha_{mat} :=
\left(\begin{matrix}\alpha_{1} & &\\ & \ddots & \\ & & \alpha_{n}
\end{matrix} \right)\, .
\]
It follows that there are formulas for $A_{i},C_{i},U,V$ in terms
of $\Delta, N^{(0)}_{mat}, N^{(k)}_{mat}, \beta^{(k)}_{mat}$,
$\alpha_{mat}$ and $B^{(1)}$. Define
\[
\tilde{\Delta} := \left(\begin{matrix} \prod_{k=1}^{r}
t_{k}^{\beta_{1}^{(k)}} & &\\ & \ddots & \\ & & \prod_{k=1}^{r}
t_{k}^{\beta_{n}^{(k)}}\end{matrix} \right)\, .
\]
We have for example
\begin{align}
A_{i} &= \beta^{(i)}_{mat} + \frac{-1}{2 \pi i} \Delta^{-1}\,
N^{(i,0)}_{mat} \Delta  \quad i \in \lbrace 1 , \ldots , r \rbrace, \notag \\
C_{i} &= [A_{i},B^{(1)}] + \frac{-1}{2 \pi i} \tilde{\Delta}^{-1}\,
N^{(i,1)}_{mat} \tilde{\Delta} + t_{i} \frac{\p}{\p t_{i}}
B^{(1)}
\quad i \in \lbrace 1 , \ldots , r \rbrace, \notag \\
C_{j} &= \frac{\p B^{(1)}}{\p t_{j}} \quad \quad j \in \lbrace r+1 ,
\ldots , m \rbrace, \notag \\
F_{\alpha} &= \frac{\p B^{(1)}}{\p y_{\alpha}} \quad \alpha \in
\lbrace 1 , \ldots , l \rbrace. \notag
\end{align}
But we will not use this formulas in the following.
Recall that an elementary section is utterly dependent on the chosen
coordinates cf. \cite{He1}. Assume that we have given two unfoldings
\begin{align}
&(\tilde{H} \rightarrow \mathbb{P}^{1} \times (M \times \mc^{l},0))
\;\; \text{with coordinates} \;\; y_{\alpha}\;\; \text{on}\;\;
\mc^{l}, \notag
\\
&(\tilde{H}' \rightarrow \mathbb{P}^{1} \times (M \times
\mc^{l'},0)) \;\; \text{with coordinates} \;\; y_{\alpha}' \;\;
\text{on}\;\; \mc^{l'}. \notag
\end{align}
We now want to examine the properties of a map
\begin{align}
&\varphi: (M \times \mc^{l'}, D \times \mc^{l'}, 0) \rightarrow ( M
\times \mc^{l} , D \times \mc^{l}, 0) \notag \\
&\text{with}\;\; (\tilde{H}', \tilde{\nabla}') \simeq
\varphi^{\ast}(\tilde{H}, \tilde{\nabla}) \notag \\
&\text{and}\;\; \varphi\mid_{M \times \lbrace 0 \rbrace} = id.
\notag
\end{align}
Therefore $\varphi$ has the following form
\begin{align}
\varphi &= (\varphi_{1}, \ldots ,\varphi_{r}, \varphi_{r+1}, \ldots
, \varphi_{m}, \varphi_{m+1}, \ldots , \varphi_{m+l}) \notag \\
&= ( t_{1} \cdot e^{u_{1}(t,y')}, \ldots , t_{r} \cdot
e^{u_{r}(t,y')}, \varphi_{r+1}, \ldots, \varphi_{m}, \varphi_{m+1},
\ldots , \varphi_{m+l}). \label{1:phiform}
\end{align}
But this has the following effect on the elementary sections (here
we identify $\varphi^{\ast}(A_{j})$ with $A'_{j}$)
\begin{align}
\varphi^{\ast} s_{j} &= \left( z^{\alpha_{j} - \frac{N^{(0)}}{2 \pi
i}} \prod_{k=1}^{r} t_{k}^{\beta_{j}^{(k)}- \frac{N^{(k)}}{2 \pi i}}
\cdot \prod_{k=1}^{r} e^{(\beta_{j}^{(k)}-\frac{N^{(k)}}{2 \pi
i})u_{k}(t,y')}\right)A'_{j} \notag \\
&= \left( \prod_{k=1}^{r} e^{(\beta_{j}^{(k)} - \frac{N^{(k)}}{2 \pi
i}) \cdot u_{k}(t,y')} \right) s'_{j} \notag\\
&= \prod_{k=1}^{r} e^{\beta_{j}^{(k)} \cdot u_{k}} \cdot \left(
\prod_{k=1}^{r} e^{- \frac{N^{(k)}}{2 \pi i} \cdot u_{k}} \right)
s'_{j}\, , \notag
\end{align}
which is in matrix notation
\begin{align}
\varphi^{\ast}(s_{1}, \ldots , s_{n}) =  \left[ \left(
\prod_{k=1}^{r} e^{- \frac{N^{(k)}}{2 \pi i} \cdot u_{k}} \right)
(s_{1}', \ldots , s_{n}') \right] \cdot \prod_{k=1}^{r}
e^{\beta_{mat}^{(k)} \cdot u_{k}} \notag \\
= (s_{1}', \ldots , s_{n}') \cdot \left(\prod_{k=1}^{r} e^{(- \frac{1}{2 \pi i} \Delta^{-1}\,
N^{(k)}_{mat} \Delta ) u_{k} }\right) \cdot \left( \prod_{k=1}^{r} e^{\beta_{mat}^{(k)} \cdot u_{k}}\right). \notag
\end{align}
We now want to expand the matrix
\[
e^{(-\frac{1}{2 \pi i} \Delta^{-1}\, \cdot N^{(k)}_{mat} \cdot
\Delta)u_{k}}
\]
in powers of $z^{-1}$. We use the matrices $N^{(k,l)}_{mat}$ defined
above. This yields
\[
\Delta^{-1} \cdot N^{(k)}_{mat} \cdot \Delta = \sum_{l \geq 0}
\tilde{\Delta}^{-1}\, N^{(k,l)}_{mat} \cdot \tilde{\Delta} \cdot
z^{-l}.
\]
Observe that for fixed $k$ the matrices $N^{(k,l)}_{mat}$  in
general do not commute. Therefore we have the following expansion:
\begin{align}
&\quad\; e^{(- \frac{1}{2 \pi i} \Delta^{-1}\, \cdot N^{(k)}_{mat} \cdot
\Delta) u_{k}} \notag\\
&= z^{0} \cdot \tilde{\Delta}^{-1}\, \cdot e^{-
\frac{1}{2 \pi i} N^{(k,0)}_{mat} \cdot u_{k}} \cdot
\tilde{\Delta} \notag \\
&+ z^{-1} \Big[ \; \frac{-1}{2 \pi i} \tilde{\Delta}^{-1}\, \cdot
N^{(k,1)}_{mat} \cdot \tilde{\Delta} \cdot u_{k} \notag \\
&\hspace{10pt}+ \frac{1}{2!} \left( \frac{-1}{2 \pi i}\right)^{2}\tilde{\Delta}^{-1}\, \cdot \left(
N^{(k,0)}_{mat} \cdot N^{(k,1)}_{mat} + N^{(k,1)}_{mat}
N^{(k,0)}_{mat} \right)
\cdot \tilde{\Delta} + \ldots \Big] \notag \\
&+ z^{-2} \cdot [\ldots] + z^{-3} \cdot [ \ldots ] +
\ldots \;. \notag
\end{align}
Now set
\begin{align}
\Gamma_{1} &= \prod^{r}_{k=1} \tilde{\Delta}^{-1}\, \cdot e^{- \frac{1}{ 2 \pi i}
N^{(k,0)}_{mat} \cdot u_{k}} \cdot \tilde{\Delta} \; \;
\text{and}
\notag \\
\Lambda &= \prod_{k=1}^{r} e^{\beta^{(k)}_{mat} \cdot u_{k}}.
\notag
\end{align}
Let $\tilde{v}_{1},\ldots, \tilde{v}_{n}$ respectively $\tilde{v}_{1}',
\ldots, \tilde{v}_{n}'$ be the canonical extensions of $v_{1},
\ldots , v_{n}$. We then have
\[
(\tilde{v}_{1}, \ldots ,\tilde{v}_{n}) =
(s_{1},
\ldots, s_{n}) \cdot \left(\sum_{k=0}^{\infty}z^{-k} \tilde{B}^{(k)}\right)  \;\; \text{with} \;\; \tilde{B}^{(0)} = I
\]
and analogously for $\tilde{H}' \rightarrow \mathbb{P}^{1} \times (
M \times \mc^{l'},0)$. Now we are able to examine how the global sections $(\tilde{v}_1, \ldots , \tilde{v}_n)$ behave under a pullback (where we identify the
flat sections $A'_{j}$ of $\tilde{H}'$ and the sections
$\varphi^{\ast}(A_{j})$ of $\varphi^{\ast}(\tilde{H})$).
\begin{align}
&\phantom{=}\;\,\,\varphi^{\ast}(\tilde{v}_{1}, \ldots , \tilde{v}_{n}) \notag \\
&=\varphi^{\ast}\left((s_{1},\ldots , s_{n}) \cdot \sum_{k=0}^{\infty}
z^{-k} \tilde{B}^{(k)} \right) \notag \\
&= (s_{1}', \ldots , s_{n}')\hspace{-1pt} \cdot \hspace{-1pt} \left(\prod_{k=1}^{r} e^{(- \frac{1}{2
\pi i} \Delta^{-1}\, N^{(k)}_{mat} \Delta) u_{k} }\right)\hspace{-2pt} \cdot \hspace{-2pt} \left(\prod_{k=1}^{r} e^{\beta^{(k)}_{mat} \cdot u_{k}} \right)\hspace{-2pt} \cdot \hspace{-2pt} \left(\sum_{k=0}^{\infty} z^{-k} (\tilde{B}^{(k)} \circ
\varphi)\hspace{-1pt}\right) \notag \\
&=:  (s_{1}', \ldots , s_{n}')  \cdot \left(z^{0} \cdot \Gamma_{1}
+ z^{-1} \cdot \Gamma_{2} + z^{-2} \cdot [ \ldots ]+ \ldots \right) \cdot \Lambda
\cdot \left( \sum_{k=0}^{\infty} z^{-k} (\tilde{B}^{(k)} \circ
\varphi) \right) \notag \\
&= (s_{1}', \ldots , s_{n}') \cdot \left[ \Gamma_{1} \cdot \Lambda + z^{-1} \cdot \left(
\Gamma_{2} \cdot \Lambda + \Gamma_{1} \cdot \Lambda \cdot (\tilde{B}^{(1)} \circ \varphi)  \right) + z^{-2} \cdot [ \ldots ] + \ldots \right] , \notag
\end{align}
and therefore
\begin{align}
& \varphi^{\ast}(\tilde{v}_{1}, \ldots ,\tilde{v}_{n}) \cdot \Lambda^{-1} \cdot
\Gamma_{1}^{-1} = \notag \\
& (s_{1}', \ldots ,
s_{n}') \cdot \left[ \mathbf{1} + z^{-1} \cdot \left( \Gamma_{2} \cdot
\Gamma_{1}^{-1} + \Gamma_{1} \cdot \Lambda \cdot  (\tilde{B}^{1}
\circ \varphi) \cdot \Lambda^{-1} \cdot \Gamma_{1}^{-1} \right) + \ldots \right] . \notag
\end{align}
Now the matrix $\Gamma_{2}$ above is defined implicitly and its form is extremely complicated.
But we do not need its explicit form in the following.\\
\noindent Recall the definition of the global sections $\tilde{v}_1, \ldots, \tilde{v}_n$, which are a canonical extension of an adapted basis $v_1, \ldots , v_n$. They were defined being elementary sections of the residue connection $\tilde{\nabla}^{res}$ on $\tilde{H}\hspace{-3pt}\mid_{\{\infty\} \times (M \times \mc^l,0)}$ when restricted to $\{\infty \} \times (M \times \mc^{l},0)$ and $(\tilde{v}_1, \ldots, \tilde{v}_n)\hspace{-3pt}\mid_{\mathbb{P}^{1} \times (M \times \{ 0 \},0)} = (v_1, \ldots, v_n)$. Now recall the correspondence between logarithmic Frobenius type structures and $(logD-trTLEP(w))$-structures. Given a $(logD-trTLEP(w))$-structure we set $K:= H\mid_{\{0 \} \times M}$ and used the canonical isomorphism between $H\mid_{\{ \infty \} \times M}$ and $K$ to shift the residue connection of $\nabla$ on $\{\infty \} \times M$ to a connection on $K$.
We therefore have
\begin{align}
(\tilde{v}_{1}, \ldots \tilde{v}_{n})\mid_{\{\infty\} \times M \times \mc^{l}} &= (s_1, \ldots , s_n)\mid_{\{\infty \} \times M \times \mc^{l}} \notag \\
(\tilde{v}_{1}', \ldots \tilde{v}_{n}')\mid_{\{\infty\} \times M \times \mc^{l'}} &= (s'_1, \ldots , s'_n)\mid_{\{\infty \} \times M \times \mc^{l'}} \notag
\end{align}
but this yields
\begin{equation}\label{1:impform}
\varphi^{\ast}(\tilde{v}_1, \ldots , \tilde{v}_n) \cdot \Lambda^{-1} \cdot \Gamma_{1}^{-1} = (\tilde{v}'_1, \ldots , \tilde{v}'_n)\; .
\end{equation}
This formula is the crucial ingredient in proving the universality of an unfolding.\\[15pt]
\noindent(\textbf{C}) We now prove the theorem \ref{U:General}.
\begin{proof}
Let $((M,0),D,K, \dr,\calc, \calu, \calv, g)$ be a Frobenius type
structure with logarithmic pole along $D$ as in theorem \ref{U:General}, let $w \in \mathbb{Z}$
and $((M,0),D, H, \nabla, P)$ be the corresponding $(logD -
trTLEP(w))$-structure. Because of the lemma and condition (IC) an
unfolding $((M \times \mc^{l},0),\tilde{D}, \tilde{H},
\tilde{\nabla}, \tilde{P})$ exists such that for the corresponding
Frobenius type structure $((M \times \mc^{l},0),\tilde{D},\tilde{K},
\tilde{\nabla}^{r},\tilde{\calc}, \tilde{\calu}, \tilde{\calv},
\tilde{g})$
\[
\tilde{\calc}_{\bullet} \xi :Der_{M \times \mc^{l}}(log
\tilde{D})_{0} \longrightarrow \tilde{K}_{0}
\]
is an isomorphism. Now we are in the isomorphism case which shows
that the constructed  Frobenius type structure corresponds to a Frobenius
manifold.\\[7pt]
It rests to prove uniqueness of this Frobenius manifold. As mentioned above
this is equivalent to the universality of the unfolding.
Therefore we now prove universality of the unfolding. Consider a second
unfolding $((M \times \mc^{l'},0),\tilde{D'},\tilde{H'},
\tilde{\nabla}',$ $\tilde{P}')$. The two connection matrices  of the
corresponding $(log\tilde{D}^{(')}-trTLEP(w))$-structures are
\begin{align}
\Omega &= \sum_{i=1}^{r} A_{i} \frac{dt_{i}}{t_{i}} + \frac{1}{z}
\sum^{r}_{i=1} C_{i} \frac{dt_{i}}{t_{i}} + \frac{1}{z} \sum_{j =
r+1}^{m} \hspace{-1pt} C_{j} dt_{j} + \frac{1}{z} \sum_{\alpha =1}^{l} F_{\alpha}
dy_{\alpha} + (\frac{1}{z^{2}} U + \frac{1}{z} V)
dz\, , \notag \\
\Omega^{'}\hspace{-1pt} &=\hspace{-1pt}  \sum_{i=1}^{r} A_{i} \frac{dt_{i}}{t_{i}} +
\frac{1}{z}\hspace{-1pt} \sum^{r}_{i=1}\hspace{-1pt} C_{i}^{'} \frac{dt_{i}}{t_{i}} +
\frac{1}{z} \hspace{-1pt} \sum_{j = r+1}^{m}\hspace{-3pt} C_{j}^{'} dt_{j} + \frac{1}{z}
\sum_{\alpha =1}^{l'} F_{\alpha}^{'} dy'_{\alpha} + (\frac{1}{z^{2}}
U^{'} + \frac{1}{z} V^{'}) dz\;  \notag
\end{align}
where $\Omega$ is with respect to adapted sections $(\tilde{v}_1, \ldots , \tilde{v}_n)$ and
$\Omega'$ is with respect to adapted sections $(\tilde{v}'_{1}, \ldots , \tilde{v}'_n)$, which are both
canonical extensions of adapted sections $(v_1, \ldots , v_n)$. \\[7pt]

\noindent We want to find a map $\varphi : (M \times \mc^{l'}, \tilde{D}') \rightarrow (M \times \mc^{l}, \tilde{D})$ with $\varphi_{\mid M \times \{ 0\} } =id$ such that
\begin{equation}\label{1:unfoldcond}
\varphi^{\ast}(\tilde{\nabla}) = \tilde{\nabla}'
\end{equation}
Now the discussion in (\textbf{B}) shows that the $\varphi^{\ast}(\tilde{v}_1, \ldots , \tilde{v}_n)\mid_{\{ \infty \} \times M \times \mc^{l'}}$ are not elementary sections of the residue connection of $\varphi^{\ast}(\tilde{\nabla})$ on $\{\infty \} \times M \times \mc^{l'}$ in general.
But we know from formula (\ref{1:impform}) that
\[
\varphi^{\ast}(\tilde{v}_1, \ldots , \tilde{v}_n) \cdot \Lambda^{-1} \cdot \Gamma_{1}^{-1} \mid_{\{ \infty \} \times M \times \mc^{l'}}
\]
are elementary sections of the residue connection of $\varphi^{\ast}(\tilde{\nabla})$. Now set
\[
\Gamma := \Lambda^{-1} \cdot \Gamma_{1}^{-1} .
\]

\noindent Therefore if we identify $(\varphi^{\ast}(\tilde{v}_{1}), \ldots ,
\varphi^{\ast}(\tilde{v}_{n})) \cdot \Gamma$ with $(\tilde{v}_{1}', \ldots , \tilde{v}_{n}')$
the following relation must be satisfied if we want (\ref{1:unfoldcond}):
\begin{equation}\label{1:conmatcond}
\Omega' = \Gamma^{-1}d \Gamma + \Gamma^{-1} \varphi^{\ast} (\Omega)
\Gamma\, .
\end{equation}
$ $\\
The reason that we have to use this matrix $\Gamma$ is, that the
sections $(\tilde{v}_{1}, \ldots, \tilde{v}_{n})$ were defined using elementary
sections (with respect to the residue connection). But elementary
sections are coordinate dependent, therefore we can not expect equality
of $(\varphi^{\ast}(\tilde{v}_{1}), \ldots , \varphi^{\ast}(\tilde{v}_{n}))$ and
$(\tilde{v}_{1}', \ldots , \tilde{v}_{n}')$, unless the first $i$
components of $\varphi$ are the identity. \\

\noindent We want to use the uniqueness statement in lemma \ref{1:longlem} that an unfolding
is determined by the first columns of the $F_{\alpha}$ in the connection matrix.
Therefore if we want to have (\ref{1:unfoldcond}), $\varphi : (M \times \mc^{l'} , \tilde{D}') \rightarrow (M \times \mc^{l}, \tilde{D})$ must satisfy $\varphi\mid_{M \times \{ 0 \} } = id$ and the first columns of the $F_{\alpha}$
must be equal.\\[7pt]
Now observe that if we have $\varphi\mid_{M \times \{ 0 \}} = id$ then we have
\[
\Gamma^{-1} d\Gamma + \Gamma^{-1}\varphi^{\ast} \left( \sum_{i=1}^{r} A_i \frac{dt_i}{t_i} \right) \Gamma \mid_{\{\infty\} \times M \times \mc^{l'}} \; = \sum_{i=1}^{r} A_i \frac{dt_i}{t_i} \mid_{\{ \infty \} \times M \times \mc^{l'}}
\]
because $(\tilde{v}'_1, \ldots , \tilde{v}'_n )$ and $\varphi^{\ast}(\tilde{v}_1, \ldots , \tilde{v}_n) \cdot \Gamma$ are both elementary sections on $\{ \infty \} \times M \times \mc^{l'}$ which coincide on $\{ \infty \} \times M \times \{ 0 \}$ and therefore coincide on all of $\{ \infty \} \times M \times \mc^{l'}$.\\[5pt]

\noindent We have $l'$ equations coming from (\ref{1:conmatcond}) and using (\ref{1:phiform})
\[
F'_{\alpha} = \sum_{i=1}^{r} \Gamma^{-1} (C_{i} \circ \varphi)
\Gamma \frac{\p u_{i}}{\p y_{\alpha}'} +\sum_{j=r+1}^{m} \Gamma^{-1}
(C_{j} \circ \varphi) \Gamma \frac{\p \varphi_{j}}{\p y_{\alpha}'} +
\sum^{l}_{\beta=1} \Gamma^{-1} (F_{\beta} \circ \varphi) \Gamma
\frac{\p \varphi_{m + \beta}}{\p y _{\beta}'} \, .
\]
We consider the first columns of $\Gamma^{-1}
\cdot C_{i}\cdot  \Gamma$ and $\Gamma^{-1} \cdot F_{\alpha} \cdot
\Gamma$ respectively and denote these with an upper
$I$. The matrices
\[
S = (C_1^I, \ldots , C_m^I,F_1^I, \ldots , F_l^I)\, , \notag \\
\]
and $\tilde{S}$ which is defined by
\[
((\Gamma^{-1} \cdot (C_{1}\circ \varphi) \cdot \Gamma)^{I} \ldots
(\Gamma^{-1} \cdot (C_{m} \circ \varphi) \cdot \Gamma)^{I}\,(\Gamma^{-1} \cdot (F_{1} \circ \varphi)
\cdot \Gamma)^{I} \ldots  (\Gamma^{-1} \cdot (F_{l} \circ \varphi) \cdot
\Gamma)^{I})
\]
are locally invertible around $\{ 0 \} \times \{ 0 \}$ because
$\Gamma \mid_{M \times \{ 0 \}} = I$, therefore $\tilde{S}\hspace{-1pt}\mid_{M
\times \{ 0 \}} = S\mid_{M \times \{ 0 \}}$ and $S$ is invertible in
$\{ 0 \} \times \{ 0 \}$ by construction. We get the following
system of equations
\[
\tilde{S}^{-1} (F'_{\alpha})^{I} = (\frac{\partial u_{1}}{\partial
y_{\alpha}'}, \ldots, \frac{\partial u_{r}}{\partial y_{\alpha}'},
\frac{\partial \varphi_{r+1}}{\partial y_{\alpha}'}, \ldots ,
\frac{\partial \varphi_{m}}{\partial y_{\alpha}'}, \frac{\partial
\varphi_{m+1}}{\partial y_{\alpha}'}, \ldots , \frac{\partial
\varphi_{m+l}}{\partial y_{\alpha}'})^{tr}\, .
\]
$ $\\[5pt]
Now $\varphi$ is determined on $M \times \{ 0 \}$ (the
identity) and we can solve for $\varphi$ with the
theorem of Cauchy-Kovalevski by using induction in $\alpha$.\\[5pt]
The beginning of the induction is clear:
\[
(\Gamma^{-1} d \Gamma + \Gamma^{-1} \varphi^{\ast} (\Omega) \Gamma)\mid_{M \times \{ 0 \}} = \Omega'\mid_{M \times \{ 0 \}}
\]
because $\Gamma\mid_{M \times \{ 0\}} = \mathbf{1}$ and $\varphi\mid_{M \times \{ 0 \}} =id$.\\
Now assume that we have constructed $\varphi$ such that
\[
(\Gamma^{-1} d \Gamma + \Gamma^{-1} \varphi^{\ast}(\Omega) \Gamma)\mid_{M \times \{ y_{\alpha} = \ldots = y_{l'} = 0 \}} = \Omega'\mid_{M \times \{ y_{\alpha} = \ldots = y_{l'} = 0 \}}.
\]
We use the equation
\[
\tilde{S}^{-1} (F'_{\alpha})^{I}\hspace{-2pt}\mid_{y_{\alpha +1}= \ldots = y_{l'} = 0} \; = (\frac{\partial u_{1}}{\partial
y_{\alpha}'}, \ldots, \frac{\partial u_{r}}{\partial y_{\alpha}'},
\frac{\partial \varphi_{r+1}}{\partial y_{\alpha}'}, \ldots ,
, \frac{\partial
\varphi_{m+l}}{\partial y_{\alpha}'})^{tr}\hspace{-1pt}\mid_{y_{\alpha +1}= \ldots = y_{l'} = 0}\,
\]
to construct an extension of $\varphi$ from $M \times \{y_{\alpha}= \ldots = y_{l'} = 0\}$ to
$M \times \{y_{\alpha +1}= \ldots = y_{l'} = 0\}$ by using the theorem of Cauchy-Kovalevski.
\noindent From the construction of the extension of $\varphi$ follows that the first column of
$F_{\alpha}'$ and the first column of the corresponding matrix in
$\Gamma^{-1} d\Gamma + \Gamma^{-1} \varphi^{\ast}(\Omega) \Gamma$  coincide on
$M \times \{y_{\alpha +1}= \ldots = y_{l'} = 0\}$. The uniqueness statement of lemma \ref{1:longlem}
in the case $l=1$ gives
\[
(\Gamma^{-1} d \Gamma + \Gamma^{-1} \varphi^{\ast}(\Omega) \Gamma)\mid_{M \times \{ y_{\alpha+ 1} = \ldots = y_{l'} = 0 \}} = \Omega'\mid_{M \times \{ y_{\alpha +1 } = \ldots = y_{l'} = 0 \}}.
\]
which was to be shown.\\[5pt] It rests to prove uniqueness of $\varphi$.
But the condition $\varphi\mid_{M \times \{ 0 \}} = id$ and the equations
\[
\tilde{S}^{-1} (F'_{\alpha})^{I} = (\frac{\partial u_{1}}{\partial
y_{\alpha}'}, \ldots, \frac{\partial u_{r}}{\partial y_{\alpha}'},
\frac{\partial \varphi_{r+1}}{\partial y_{\alpha}'}, \ldots ,
\frac{\partial \varphi_{m}}{\partial y_{\alpha}'}, \frac{\partial
\varphi_{m+1}}{\partial y_{\alpha}'}, \ldots , \frac{\partial
\varphi_{m+l}}{\partial y_{\alpha}'})^{tr}\, .
\]
show that $\varphi$ is unique.\\
\end{proof}

\section{Logarithmic Frobenius manifolds and quantum
cohomology}\label{S: qcoh} $ $\\
In this section we show how to
construct a Frobenius manifold out of quantum cohomology where we
are following \cite{Ma}. With the definition of a logarithmic
Frobenius manifold at hand we show how this construction extends for
smooth, projective varieties to give us
logarithmic Frobenius manifolds. This construction enables us to
prove a partial generalization of the first Reconstruction theorem
in \cite{KM}.
\subsection{Quantum Cohomology and Frobenius manifolds} $
$\\[10pt]
We now define quantum cohomology on the even dimensional cohomology
ring. That means we define a potential on $H^{*}_{even}(X,\mc)$
which we use to define a quantum deformation of the cup product. Let
$T_{0}=1 \in H^{0}(X,\mc)$, $T_{1},\ldots,T_{r}$ be a basis of
$H^{2}(X,\mc)$ and let $T_{r+1},\ldots, T_{m}$ be a basis for the
other cohomology groups lying in $H^{\ast}_{even}(X,\mc)$. Put $ \gamma = \sum_{i=0}^{m} t_{i}
T_{i}$. In this section we only consider the \textbf{free} parts of $H_2(X,\mz)$ and $H^2(X,\mz)$ and denote them by the same letters.\\[10pt]
\begin{definition}
Let $X$ be a smooth projective variety. Then the Gromov-Witten
potential is the formal sum
\begin{align}
\Phi(\gamma) &= \sum_{n=0}^{\infty} \sum _{\beta \in H_{2}(X,\mz)}
\frac{1}{n!} \langle I_{0,n,\beta} \rangle (\gamma^{n})
\label{2:infsum2} \\
&= \frac{1}{6} \int_{X} \gamma^{3} + \sum_{n=0}^{\infty} \sum
_{\substack {\beta \in H_{2}(X,\mz)\\ \beta \neq 0}} \frac{1}{n!}
\langle I_{0,n,\beta} \rangle (\gamma^{n})  \notag \\
&= \; \Phi_{class} \hspace{7pt} + \; \Phi_{quantum} \notag
\end{align}
where we set $\langle I_{0,n,\beta} \rangle = 0$ for $n\leq 2$.
\end{definition}
\noindent There are two problems with this definition. The first is, if for fixed $n$ the sum
\begin{equation}\label{2:infsum1}
\sum_{\beta \in H_2(X,\mz)}\frac{1}{n!} \langle I_{0,n,\beta} \rangle (\gamma^n)
\end{equation}
does converge. Usually one introduces the Novikov ring as the coefficient ring to split the contribution of the different $\beta$. But here we will assume that $(\ref{2:infsum1})$ does converge.
The second problem is if the whole sum (\ref{2:infsum2}) gives a well-defined function at least in some domain of $H^{\ast}(X,\mc)$.  \\
\noindent \textbf{Assumption:} We assume in the rest of the paper
that the individual summands in the Gromov-Witten potential are convergent. We assume also that the Gromov-Witten potential gives a well defined holomorphic function in the domain
\[
B:= B_0 \times \prod_{i=1}^r \{Re(t_i) < r_i \} \times \prod_{j=r+1}^m B_j \subset H^{\ast}_{even}(X,\mc) \; ,
\]
where $B_0$ resp. $\prod_{j=r+1}^m B_j$ are (poly)discs in $H^0$ resp. $H^4 \oplus \ldots \oplus H^{2 dim X}$ and $r_i \in \mathbb{R}$ for all $i$.\\[10pt]

\begin{remark}\label{2:Potdiv}
Write $\gamma = \gamma_{0} + \delta$ with $\delta \in H^{2}(X,
\mc)$. We use the divisor axiom and the fundamental class axiom to rewrite the expression for the
potential:
\begin{align}
\Phi_{quantum}(\gamma_{0} + \delta) &= \sum_{i,k \geq 0} \sum_{\beta
\neq 0} \frac{1}{i! k!} \langle I_{0,n,\beta} \rangle
(\gamma_{0}^{\otimes i} \otimes \delta^{\otimes k})\notag \\
&= \sum_{i \geq 0} \sum_{\beta \neq 0} \frac{e^{(\beta,
\delta)}}{i!} \langle I_{0,i,\beta} \rangle (\gamma_{0}^{\otimes
i})\notag \\
&=\sum_{j_{r+1},\ldots , j_{m}} \sum_{\beta \neq 0} e^{(\beta,
\delta)} \langle I_{0,i,\beta} \rangle (
T_{r+1}^{j_{r+1}}\otimes \ldots \otimes
T_{m}^{j_{m}})\frac{t_{r+1}^{j_{r+1}} \ldots
t_{m}^{j_{m}}}{j_{r+1}! \ldots j_{m}!}. \notag
\end{align}
\end{remark}
$ $\\

\begin{definition}\label{2:defnmult}
Let $\Phi$ be the Gromov-Witten potential for a smooth projective
variety X. Then define
\[
T_{i} \ast T_{j} = \sum_{k} \frac{\partial^{3}\Phi}{\partial t_{i}
\partial t_{j} \partial t_{k}} T^{k}
\]
where $\{ T^k \}$ is the dual base with respect to the Poincar\'{e} pairing, i.e
$\int_{X} T_i \cup T^j = \delta_{ij}$.
Extending this linearly gives the big quantum product on the
cohomology $H^{*}_{even}(X,\mc)$.
\end{definition}
\begin{lemma}\label{2:part3phi}
For all $i,j,k,$ we have
\begin{align}
&\frac{\partial^{3} \Phi}{\p t_{i} \p t_{j} \p t_{k}} =
\sum_{n=0}^{\infty} \sum_{\beta \in H_{2}(X, \mz)} \frac{1}{n!}
\langle I_{0,n+3, \beta} \rangle (T_{i},T_{j},T_{k},
\gamma^{n}) \notag \\
&=\sum_{j_{r+1},\ldots , j_{m}} \sum_{\beta \neq 0} e^{(\beta,
\delta)} \langle I_{0,i+3,\beta} \rangle (T_{i} \otimes T_{j}
\otimes T_{k} \otimes T_{r+1}^{j_{r+1}}\otimes
\ldots \otimes T_{m}^{j_{m}})\frac{t_{r+1}^{j_{r+1}}
\ldots t_{m}^{j_{m}}}{j_{r+1}! \ldots j_{m}!} \notag \\
& \qquad +  \int_{X} T_{i} \cup T_{j} \cup T_{k}\; .
\notag
\end{align}
\end{lemma}
\begin{proof}
See \cite{CK} lemma 8.2.3 .
\end{proof}
$ $\\
\noindent It is a well known fact that $\ast$ is commutative and
associative, with unit $T_{0}$. For a proof see \cite{CK}.\\
$ $\\
\subsubsection{Construction of Frobenius manifold} $
$\normalsize \\[10pt]
Because of the assumption on the Gromov-Witten potential stated above we will
get a Frobenius manifold on the open subset $B \subset H^{\ast}_{even}(X,\mc)$.
We consider now $B \subset H^{*}_{even}(X,\mc)$ as a manifold with
global coordinates $t_{i}$ with respect to a basis $\lbrace T_{i}
\rbrace$ as above. At each point the tangent space is canonically
isomorphic to $H^{*}_{even}(X,\mc)$. Denote the global vector fields as
$\pti$ and define a metric
\[
g(\pti,\ptj) := g_{ij} = \int_{X} T_{i} \cup T_{j}.
\]
Observe that the metric is constant so that the induced Levi-Civita
connection is $\nabla = d$ with respect to the sections $\pti$.
Denote by $deg (\pti)$ the degree of $\pti$ in $H^{*}(X,\mc)$.
We define the Euler vector field as
\[
E := \sum_{i} ( 1 - \frac{deg (\pti)}{2} ) t_{i} \pti + \sum_{deg
(\ptj) =2 } r^{j} \ptj
\]
where the $r^{j}$ are defined as
\[
c_{1}(X) = \sum_{deg T_{j} =2 } r^{j} T_{j}\, .
\]
\begin{proposition}
Let $E$ as above. It holds
\begin{align}
&a.) E \Phi_{quantum} = (3- dim X) \Phi_{quantum}\, , \notag \\
&b.) (E - E(0)) \Phi_{class} = (3 - dim X) \Phi_{class}\, , \notag \\
&c.)\, E \Phi_{ijk} = (3 -dim X) \Phi_{ijk} - \left( 3 - \frac{deg \p_{t_i} + deg \p_{t_j} + deg \p_{t_k}}{2} \right) \Phi_{ijk} \, . \notag
\end{align}
\end{proposition}
\begin{proof}
Recall that
\[
\Phi_{quantum} = \sum_{j_{r+1},\ldots , j_{m}} \sum_{\beta
\neq 0} \frac{e^{(\beta, \delta)}}{i!} \langle I_{0,i,\beta} \rangle
(T_{r+1}^{j_{r+1}}\otimes \ldots \otimes
T_{m}^{j_{m}})\frac{t_{r+1}^{j_{r+1}} \ldots
t_{m}^{j_{m}}}{j_{r+1}! \ldots j_{m}!}.
\]
$ $\\[5pt]
We apply $E$ to that term. The $E(0)$ term acts only upon
$e^{(\beta, \delta)}$ and multiplies it by $(c_{1}(X), \beta)$. The
$E - E(0)$ part multiplies any monomial $t_{a_{1}}^{j_{1}} \ldots
t_{a_{n}}^{j_{n}}$ in non-divisorial coordinates by $\sum_{i}(1-
\frac{deg(T_{a_{i}})}{2})$. Observe that the Gromov-Witten
invariants contribute only for such $\beta$ which satisfy
\[
\frac{1}{2} \sum_{i =1}^{n} deg(T_{a_{i}}) = dim X + \int_{\beta}
c_{1}(X) + n -3.
\]
Hence every non-vanishing term is an eigenvector of $E$ with
eigenvalue $ 3- dim X$. \\$b.)$ is clear. \\
$c.)$ We are using $[E, \p_{t_m}] = - (1 - \frac{deg \p_{t_m}}{2})\p_{t_m}$:
\begin{align}
E \Phi_{ijk} &= E \p_{t_i} \p_{t_j} \p_{t_k} \Phi \notag \\
&= \p_{t_i} \p_{t_j} \p_{t_k} E \Phi - \left(3 - \frac{deg \p_{t_i} + deg \p_{t_j} + deg \p_{t_k}}{2} \right) \Phi_{ijk}\, . \notag
\end{align}
The statement follows with
\begin{align}
E \Phi &= E \Phi_{quantum} + (E-E(0))\Phi_{class} + E(0)\Phi_{class} \notag \\
&= (3 - dim X) \Phi_{quantum} + (3-dim X) \Phi_{class} + E(0)\Phi_{class} \notag
\end{align}
and the observation that $\p_{t_i} \p_{t_j} \p_{t_k} E(0) \Phi_{class} = 0$ because $\Phi_{class}$ is only cubic in the $t_m$.
\end{proof}
\noindent We now have to check that this defines the structure of a
Frobenius manifold. The properties $(1)-(4)$ in definition \ref{1:DefFrob} are easy to show. We are using Einstein sum convention.\\[8pt]
$(5)$\; $Lie_{E}(\ast)= 1 \cdot \ast$ and $ Lie_{E}{g} = (2-d) \cdot g$ with $d= dim X$.\\
Here we use proposition 2.10c .
\begin{align}
&\quad \;\, Lie_{E}(\ast)(\pti,\ptj)\notag \\
&= [E, \pti \ast \ptj] - [E, \pti] \ast
\ptj - \pti \ast [E, \ptj] \notag \\
&= [E, \Phi_{ijk} g^{lk} \ptl] + (1 - \frac{ deg(\pti)}{2}) \pti
\ast \ptj
+ (1 - \frac{ deg(\ptj)}{2}) \pti \ast \ptj \notag \\
&= (3 - dim X) \pti \ast \ptj - (1 - \frac{
deg(\ptk)}{2}) \Phi_{ijk} g^{kl} \ptl - (1 - \frac{ deg(\ptl)}{2})
\Phi_{ijk}
g^{kl} \ptl \notag \\
&= \pti \ast \ptj , \notag \\
& \notag \\
&\quad \;\, Lie_{E}(g)(\pti,\ptj)\notag \\
&= E g(\pti,\ptj) - g([E,\pti],\ptj) -
g(\pti,[E,\ptj]) \notag \\
&= 0 + (1 - \frac{deg(\pti)}{2}) g_{ij} + (1 - \frac{deg(\ptj)}{2})
g_{ij} \notag \\
&= (2 - dim X) g(\pti,\ptj) . \notag
\end{align}

\subsubsection{Logarithmic Frobenius manifolds}\label{2:LFM} $
$\normalsize \\[10pt] Above we constructed a Frobenius manifold
on an open subset $B$ of $H^{*}_{even}(X, \mc)$ but this description has a little drawback. The
point where the multiplication is simplest is not included. We
now consider smooth, projective varieties $X$.
Choose as above a homogeneous basis $T_{0} \in H^{0},\;
T_{1}, \ldots, T_{r} \in H^{2}, T_{r+1} \in H^{4}, \ldots, T_{m} \in
H^{2 dim X}$ such that the
$T_{1}, \ldots , T_{r}$ are homogeneous with respect to the decomposition $H^2 = H^{2,0} \oplus H^{1,1} \oplus H^{0,2}$. The part of the basis which lies in $H^{1,1}(X)$ shall also lie in the K\"ahler cone.\\
If we look at
\[
\Phi_{quantum}=\sum_{j_{r+1},\ldots , j_{m}} \sum_{\beta \neq
0} \frac{e^{(\beta, \delta)}}{i!} \langle I_{0,i,\beta} \rangle
(T_{r+1}^{j_{r+1}}\otimes \ldots \otimes
T_{m}^{j_{m}})\frac{t_{r+1}^{j_{r+1}} \ldots
t_{m}^{j_{m}}}{j_{r+1}! \ldots j_{m}!}\, ,
\]
we see that the quantum potential is periodic with respect to
$\delta$ but observe that $\Phi_{class}$ is not.
Consider the map
\begin{align}
\varphi : B &\mapsto M^{\ast}:=B_0 \times \prod_{i=1}^r \{ 0 < |q_i| < e^{r_i} \} \times \vspace{-4pt}\prod_{j=r+1}^m B_j\, , \notag \\
\varphi : \hspace{6 pt} t_{i} &\mapsto q_{i} = e^{t_{i}} \quad i \in \{1, \ldots ,r\}\, ,\notag \\
\varphi : \hspace{5 pt} t_j &\mapsto t_j \qquad \quad\,\, j \in \{0,r+1,\ldots ,m\}\, , \notag \\
\varphi_{*}: \pti &\mapsto q_{i} \pqi .  \notag
\end{align}
The potential $\Phi_{quantum}$ is well-defined on the image but
$\Phi_{class}$ is not. But if we only consider the third derivatives
of the potential, $\p^{3} \Phi_{class}$ is a constant and because of
that well-defined on the image.\\
We want to investigate the limit $q_{i} \rightarrow 0$ for some $i$.
Therefore we have to look at the term $e^{(\beta,\delta)}$ more
closely. Let $ \delta = \sum_{i=1}^{r} t_{i} T_{i}$, then we have
\[
e^{(\beta, \delta)} = \prod_{i}^{r} q_{i}^{(\beta,T_i)}.
\]
Now for all effective $\beta$ we have $(\beta,T_{i}) \geq 0$ but
only these $\beta$ occur in the Gromov-Witten potential. Therefore
the limit is defined.\\
In the limit $q_{i} \rightarrow 0 \; \forall i$ all terms in the
quantum potential vanish because $\beta \neq 0$ and the $T_{i}, \; i
\in \{ 1 ,\ldots ,r \}$ are a basis of $H^{2}(X)$.\\
Set
\[
M = B_0 \times \prod_{i=1}^{r}\{ 0 \leq |q_i| < 1 + \epsilon \} \times \prod_{j=r+1}^m B_j \;.
\]
It is now clear that $\Phi_{quantum}$, being a power series in the $q_i$ for $i \in \{1, \ldots ,r \}$, is convergent on $M$.
As a consequence we have defined a Frobenius manifold structure on $M^{\ast}$
and the quantum multiplication degenerates at $q_i = 0 \, \forall i$ to the usual cup product.
It rests to examine if this defines a logarithmic
Frobenius manifold.
Note that the divisor $D = M \setminus M^{\ast}$ is
normal crossing and that the $\varphi_{\ast}(\p_{t_i})$ are a basis of $\dmldc$.  Because of $
\pti \mapsto q_{i} \pqi$ we have
\[
E = t_{0} \partial_{t_{0}} + \sum_{i=r+1}^{m}(1 -
\frac{deg({\pti})}{2}) t_{i} \pti + \sum_{j=1}^{r} r^{j} q_{j} \pqj
\]
so $e, E \in \dmld$.\\
Because $\Phi_{ijk}$ is defined for $q_{i} \rightarrow 0$ and
$\varphi$ is locally biholomorphic we can define a multiplication on
$\dmld$
\[
\varphi_{*}(\pti) \ast \varphi_{\ast}(\ptj) := \varphi_{*}(\pti \ast
\ptj).
\]
Because $g(\pti,\ptj)$ is constant, we can define
\[
g(\varphi_{*}(\pti),\varphi_{*}(\ptj)):= g(\pti,\ptj).
\]
We extend both the multiplication and the metric g
$\calo_{M}$-linearly. It follows that $M$ carries the structure of a
logarithmic Frobenius manifold. We call the point $p$ defined by
$\{t_0=q_1 = \ldots = q_r = t_{r+1} = \ldots t_{m} =0 \}$ the large radius limit point.\\[5pt]

\subsection{The first Reconstruction Theorem} $ $\normalsize
\\[10pt] \noindent We want to give a geometric interpretation and
partial generalization of the first Reconstruction Theorem of
Kontsevich and Manin given in [KM]. For the convenience of the
reader we restate the theorem.\\[-5pt]
\begin{theorem}\label{3:KM}
Let $X$ be a smooth projective variety with the property that
$H^{*}(X ,\mq)$ is generated by $H^{2}(X, \mq)$. Also assume that we
know the Gromov-Witten invariants $\langle I_{0,3,\beta} \rangle
(\alpha_{1}, \alpha_{2}, \alpha_{3})$ for all $\beta \in H_{2}(X,
\mz)$ with $deg \, \alpha_{3} =2$. Then we can determine all tree-level
Gromov-Witten classes $I_{0,n,\beta}(\alpha_{1}, \ldots ,
\alpha_{n})$ for all $\beta \in H_{2}(X, \mz)$.
\end{theorem}
$ $\\
\noindent  In \cite{KM} the additional condition $\int_{\beta} c_{1}(X) \leq dim X + 1$ is demanded. This can be explained as follows. In order to be non-zero the classes must
fulfill $\sum deg(\alpha_{i}) = 2 ( \int_{\beta} c_{1}(X) + dim
X)$ but we have $\sum_{i=1}^3
deg(\alpha_{i}) \leq 4\,
dim X +2$ (assuming $deg(\alpha_{3}) =2$). Therefore this condition is automatically fulfilled.\\
There is the following partial generalization of the theorem:\\
\begin{theorem}\label{3:myKM}
Let $\hspace{-1pt}X\hspace{-2pt}$ be a smooth projective variety. Let
$\delta_{1},\ldots,\hspace{-1pt} \delta_{u}\hspace{-1pt}$ be generators of $\hspace{-1pt}H^{*}_{even}(X,\mc)$ with respect to the cup product. Set $\sum_{k=1}^u \mc \cdot \delta_{k} = W \subset
H^{*}_{even}(X,\mc)$. Assume that we know the Gromov-Witten invariants $\langle I_{0,n,\beta} \rangle
(\alpha_{1}, \alpha_{2}, \alpha_{3}, \ldots , \alpha_{n})$ for all $\beta \in H_{2}(X,
\mz)$ with $\alpha_{3}, \ldots , \alpha_{n} \in W$ and for all $n \geq 3$. Then we can determine all tree-level Gromov-Witten classes $I_{0,n,\beta} (\alpha_{1}, \ldots , \alpha_{n})$ for all $\beta \in H_2(X, \mz)$.
\end{theorem}
\begin{proof}
We first remark that the tree-level Gromov-Witten classes can be uniquely reconstructed
from the tree-level Gromov-Witten invariants. This follows from \cite{KM} proposition $2.5.2$ .
Therefore it rests to prove that we can reconstruct the corresponding Gromov-Witten invariants.
Recall the map
\[
\varphi:\; B \mapsto M^{\ast} .
\]
and the logarithmic Frobenius manifold constructed on $M$. In proposition \ref{1:FrobmfFrobstruc} we proved that the latter gives rise to a logarithmic Frobenius type structure on $M$ . Restrict this logarithmic Frobenius type structure to the submanifold $N = \overline{\varphi(W \cap B)}$.  The Higgs field on $N$ (which comes from the quantum multiplication restricted to $N$) incorporates only Gromov-Witten invariants of the form $\langle I_{0,n,\beta} \rangle
(\alpha_{1}, \ldots , \alpha_{n})$ with $\alpha_{3}, \ldots , \alpha_{n} \in W$ (cf. definition \ref{2:defnmult} and lemma \ref{2:part3phi}). Let $\xi = \varphi_{\ast}(\p_{t_0})=\p_{t_0}$. Then it is clear that $(IC)$ is satisfied. The condition $(GC)$ holds
because the $\delta_{1}, \ldots, \delta_{u}$ are generators of
$H^{*}_{even}(X,\mc)$ and $\ast$ is simply the cup product at $p$ (the large radius limit point). Recall the definition of $E$ in the quantum cohomology case and that of $\calv$. We see that
$\nabla_{\partial_{t_{0}}} E = \partial_{t_{0}}$. So $(EC)$ is satisfied. Now theorem
(\ref{U:General}) states that we have a universal unfolding which is
equivalent to a Frobenius manifold. So we have reconstructed the big
quantum product. To extract the Gromov-Witten invariants we recall
the following formula
\[
T_{i} \ast T_{j} = \sum_{k} \frac{\p^{3} \Phi}{\p t_{i} \p t_{j} \p
t_{k}} T^{k}\, .
\]
Therefore we know $\frac{\p^{3} \Phi}{\p t_{i} \p t_{j} \p t_{k}}$
for all $i, j, k$. If we integrate this three times we get $\Phi$.
Now recall the definition of $\Phi$:
\[
\Phi(\gamma) = \sum_{n=0}^{\infty} \sum_{\beta \in H_{2}(X,\mz)}
\frac{1}{n!} \langle I_{0,n,\beta} \rangle (\gamma^{n})\, .
\]
Notice that because of the integration we have an ambiguity in the
terms for $n \leq 2$ but with the definition above ($\langle
I_{0,n,\beta}\rangle = 0$ for $n\leq 2, \beta \neq 0$) this ambiguity is fixed.\\[5 pt]

\noindent Now we want to extract the single Gromov-Witten invariants out of the potential.
Recall remark \ref{2:Potdiv} where we saw that
\[
\Phi_{quantum}(\gamma_{0} + \delta) =\sum_{j_{r+1},\ldots , j_{m}} \sum_{\beta \neq 0} e^{(\beta,
\delta)} \langle I_{0,i,\beta} \rangle (
T_{r+1}^{j_{r+1}}\otimes \ldots \otimes
T_{m}^{j_{m}})\frac{t_{r+1}^{j_{r+1}} \ldots
t_{m}^{j_{m}}}{j_{r+1}! \ldots j_{m}!}.
\]
Because of our assumption that the potential is holomorphic, we can differentiate with respect to $t_i$, $j_i$ times
such that $\sum_{i=r+1}^n j_i = n$ to get the term
\[
\frac{\p^n \Phi_{quantum}}{\p^{j_{r+1}}t_{r+1} \ldots \p^{j_m}t_m}(\delta,0) = \sum_{\beta \neq 0} e^{(\beta,
\delta)} \langle I_{0,n,\beta} \rangle (
T_{r+1}^{j_{r+1}}\otimes \ldots \otimes
T_{m}^{j_{m}}).
\]
Choose a basis $T_1, \ldots T_r$ in $H^2(X,\mz)$ and a basis $S_1, \ldots S_r$ in $H_2(X,\mz)$
such that $(S_i,T_j) = \delta_{ij}$, this is possible due to Poincar\'{e} duality. Restrict the term above to the
$\mr$-vector space spanned by $i T_1, \ldots, i T_r$ and let $\mathbf{t} =(t_1, \ldots ,t_r)$ be coordinates on that space. Then the term above
gets a multi-dimensional Fourier series
\[
\sum_{\mathbf{l} \in \mz^d\setminus \{ 0 \}} e^{i\, \mathbf{l} \cdot \mathbf{t}}\langle I_{0,n,\beta} \rangle (
T_{r+1}^{j_{r+1}}\otimes \ldots \otimes
T_{m}^{j_{m}}).
\]
Now the single coefficients can be obtained by integration.
\end{proof}
$ $\\[-5pt]
Let us compare the two theorems. If $H^{\ast}_{even}(X,\mc)$ is $H^{2}$-generated, then
$W = H^{2}$. Thus the corresponding Higgs field incorporates only
Gromov-Witten invariants of type $\langle I_{0,3,\beta} \rangle
(T_i \otimes T_j \otimes T_k)$ with $deg(T_k) =2 $  and the other
restrictions mentioned above (here we use the divisor axiom). Therefore theorem (\ref{3:myKM}) generalizes
theorem (\ref{3:KM}) for smooth projective varieties which fulfill the convergence assumptions above.\\

\section{Hodge Asymptotics}\label{S: hodgeasym} $
$\\
In this section we prove the existence of a filtration $\mcu_{\bullet}$
opposite to the Hodge filtration $\mcf^{\bullet}$ of a polarized
variation of Hodge structures which degenerates along a normal crossing divisor (prop. \ref{3:oppFil}).
This is used in theorem \ref{3:VPHSFrob2} to show that the above data give rise to a logarithmic Frobenius manifold if the limiting PMHS is Hodge-Tate, split over $\mr$ and an additional generation condition is satisfied.

\subsection{Generalities} $ $ \\
 First we recall some definitions.
\begin{definition}
A mixed Hodge structure is a triple $(H_{\mz}, W_{\bullet},
F^{\bullet})$, where $H_{\mz}$ is a lattice, $W_{\bullet}$ is a
finite increasing filtration on $H_{\mq} = H_{\mz} \otimes \mq$, and
$F^{\bullet}$ is a finite decreasing filtration on $H = H_{\mz}
\otimes \mc$, such that the induced filtration on the quotient
$Gr^{W}_{k}H = W_{k}/W_{k-1}$ defines a Hodge structure of weight
$k$. Here $F^{p}Gr^{W}_{k}H$ is the image of $F^{p} \cap W_{k}$ in
$Gr^{W}_{k}H$,
\[
F^{p}Gr^{W}_{k}H = (W_{k} \cap F^{p} + W_{k-1}) / W_{k-1} \; .
\]
$F^{\bullet}$ is called the Hodge filtration and $W_{\bullet}$ the
weight filtration.
\end{definition}
Let $H $ be as above and $H_{\mr} :=  H_{\mz} \otimes \mr$. Let $S$ be a nondegenerate
$(-1)^{w}$-symmetric pairing on $H$ with real values on $H_{\mr}$.
Let $N$ be a nilpotent endomorphism of $H_{\mr}$ and of $H$ which is
an infinitesimal isometry of $S$.  One obtains a weight filtration
$W_{\bullet}$ in the following way:
\begin{lemma}
Let $(H,H_{\mr},S,N,w)$ as above.
\renewcommand{\labelenumi}{{\normalfont (\alph{enumi})}}
\begin{enumerate}
\item There exists a unique finite increasing filtration
$W_{\bullet}$ on $H_{\mr}$ such that\\ $N(W_{l}) \subset W_{l-2}$ and
such that $N^{l}:Gr^{W}_{w+l} \longrightarrow Gr^{W}_{w-l}$ is an
isomorphism.
\item The filtration satisfies $S(W_{l},W_{l'}) = 0\;$ for $l + l' < w$.
\item A nondegenerate $(-1)^{w+l}$-symmetric bilinear form $S_{l}$
is well defined on $Gr^{W}_{w+l}$ for $l \geq 0$ by $S_{l}(a,b):=
S(\tilde{a},N^{l} \tilde{b})$ for $a,b \in Gr^{W}_{w+l}$ with
representatives $\tilde{a},\tilde{b} \in W_{w+l}$.
\item The primitive subspace $P_{w+l} \subset Gr^{W}_{w+l}$ is
defined by
\[
P_{w+l} := ker(N^{l+1} : Gr^{W}_{w+l} \longrightarrow
Gr^{W}_{w-l-2})
\]
for $l \geq 0$ and by $P_{w+l} := 0$ for $l < 0$. Then
\[
Gr^{W}_{w+l} = \bigoplus_{i \geq 0} N^{i} P_{w+l+2i}
\]
and this decomposition is orthogonal with respect to $S_{l}$ if $l
\geq 0$.
\end{enumerate}
\end{lemma}
\begin{proof}
See \cite{Sch} Lemma 6.4 .
\end{proof}
\begin{definition}
A polarized mixed Hodge structure of weight $w$ consists of
$(H,H_{\mr},S,N,w)$, $W_{\bullet}$ as above and an exhaustive
decreasing Hodge filtration $F^{\bullet}$ on $H$ with the following
properties.
\begin{enumerate}
\item $F^{\bullet}Gr^{W}_{k}$ gives a pure Hodge structure of weight
$k$, that means,
\[
Gr^{W}_{k} = F^{p}Gr^{W}_{k} \oplus \overline{F^{k+1-p}Gr^{W}_{k}}.
\]
\item $N(F^{p}) \subset F^{p-1}$, i.e. $N$ is a $(-1,-1)$-morphism
of mixed Hodge structures.
\item $S(F^{p},F^{w+1-p}) = 0$.
\item The pure Hodge structure $F^{\bullet}P_{w+l}$ of weight $w+l$
on $P_{w+l}$ is polarized by $S_{l}$, that means
\[
S_{l}(F^{p}P_{w+l},F^{w+l+1-p} P_{w+l}) = 0,
\]
\[
i^{p-(w+l-p)}S_{l}(a, \overline{a}) > 0 \quad \text{for} \; a \in
F^{p} P_{w+l} \cap \overline{F^{w+l-p} P_{w+l}} - \{ 0 \}\, .
\]
\end{enumerate}
\end{definition}
\noindent A PMHS comes equipped with a generalized Hodge
decomposition, namely Deligne's $I^{p,q}$.
\begin{lemma}\label{3:LemIpq}
For a PMHS as above define
\begin{align}
I^{p,q} &:= (F^{p} \cap W_{p+q}) \cap (\overline{F^{q}} \cap W_{p+q}
+ \sum_{j > 0} \overline{F^{q-j}} \cap W_{p+q-j-1}), \notag
\\
I^{p,q}_{0} &:= ker(N^{p+q-w+1} : I^{p,q} \longrightarrow
I^{w-q-1,w-p-1}). \notag
\end{align}
Then
\begin{align}
F^{p} &= \bigoplus_{i,q: i \geq p} I^{i,q}, &W_{l} &=
\bigoplus_{p+q \leq l} I^{p,q}, \notag \\
N(I^{p,q}) &\subset I^{p-1,q-1}, &I^{p,q} &= \bigoplus_{j \geq 0}
N^{j} I_{0}^{p+j,q+j} \notag
\end{align}
\begin{align}
S(I^{p,q},I^{r,s}) = 0 \quad &\text{for} \; (r,s) \neq (w-p,w-q),
\notag \\
S(N^{i} I^{p,q}_{0}, N^{j} I^{r,s}_{0}) = 0 \quad &\text{for} \;
(r,s,i+j) \neq (q,p,p+q-w), \notag
\end{align}
\begin{align}
I^{q,p} \simeq \overline{I^{p,q}} \quad &\text{mod } W_{p+q-2},
\notag \\
I^{q,p}_{0} \simeq \overline{I^{p,q}_{0}} \quad &\text{mod }
W_{p+q-2}. \notag
\end{align}
\end{lemma}
\begin{proof}
See \cite{De1} (1.2.8), \cite{CaK} or \cite{He3} 2.3 .
\end{proof}

\begin{definition}
A variation of Hodge structures of weight $n$ is a quadruple
$(\mch_{\mz}, \nabla,$ $M , \mcf^{\bullet} )$ where $M$ is a complex
manifold, $\mch_{\mz}$ is a local system with coefficient group
$\mz^{h}$ for some $h$, and $\nabla$ a flat holomorphic connection
on the holomorphic vector bundle $\mch = \mch_{\mz} \otimes
\calo_{M}$ such that the sections of $\mch_{\mz}$ are $\nabla$-flat. $\mathcal{F}^{\bullet}$ is a decreasing filtration of
holomorphic vector bundles $\mathcal{F}^{p} \subset \mathcal{H}$.
These objects have to satisfy the following two conditions:
\begin{itemize}
\item for each point $t \in M$ the filtration
$\mathcal{F}^{\bullet}$ induces a filtration
$\mathcal{F}^{\bullet}_{t}$ on the fibre $\mch_{\mz} \otimes \mc$ at
$t \in M$ constituting an HS of weight $n$,
\item $\nabla (\mathcal{F}^{p}) \subset \Omega^{1}_{M}
\otimes_{\calo_{M}} \mathcal{F}^{p-1}$ for each $p$.
\end{itemize}
Moreover, if there is a flat, non-degenerate bilinear form $S:
\mch_{\mz} \times \mch_{\mz} \longrightarrow \mq_{M}$ defining a
polarized HS for every $t \in M$, then the VHS $(\mch_{\mz}, \nabla,
M , \mcf^{\bullet} )$ is said to be polarized.
\end{definition}
$ $\\
Now assume that we have given a manifold $M = \Delta^{n}$, a normal crossing
divisor $Z$ on $M$  such that $ M \setminus Z = (\Delta^{\ast})^{r} \times \Delta^{n-r}, Z = \bigcup^{r}_{j=1} D_j$ and a polarized VHS on $M \setminus Z$. We want
to construct a Frobenius manifold which encodes these data. First we
construct a filtration $\mcu_{\bullet}$ which is opposite to the
Hodge filtration $\mcf^{\bullet}$. This is done using Schmid's
limiting MHS and Deligne's $I^{p,q}$'s. Then following \cite{HM}
chapter 5 we construct a Frobenius manifold out of the latter.\\
Now let $ H \rightarrow M \setminus Z$ be the flat vector bundle with connection $\nabla$ coming
from a polarized VHS. Let $\mct_j$ be the mondromies of $\nabla$ with respect to $D_j$ for $j = 1, \ldots , r$. Let $N_j := \log(\mct_{j,u})$ be the logarithm of the unipotent part of $\mct_j$. Observe that both $\mct_j$ and $N_j$ act on the fibers of $H$ and all $\mct_i$ and $N_j$ commute with each other. Let
\begin{align}
\tau : U &\longrightarrow M \setminus Z \notag \\
(z,w) = (z_1, \ldots , z_r, w_{r+1}, \ldots , w_n) &\mapsto (e^{2 \pi i z_1}, \ldots , e^{2 \pi i z_r}, w_{r+1}, \ldots , w_n) = (t,w) \notag
\end{align}
be the universal covering. We fix a flat isomorphism for the rest of the paper
\[
\rho_{\infty}: \tau^{\ast}H \simeq H^{\infty} \times U \, .
\]
For each $(z,w) \in U$ we get an isomorphism
\[
\rho_{(z,w)} : H_{\tau(z,w)} \longrightarrow H^{\infty} \, .
\]
Now we have
\[
\rho_{(z_1, \ldots , z_j + 1 ,\ldots , z_r , w)} \circ \mct_j = \rho_{(z,w)} \, .
\]
The monodromies $\mct_j$ act on $H^{\infty}$ through the maps
\[
\rho_{(z,w)} \circ \mct_j \circ (\rho_{(z,w)})^{-1} : H^{\infty} \rightarrow H^{\infty}
\]
which are actually independent of $(z,w)$. Therefore we denote them by $\mct_j$, too.
The same applies to the $N_j$.\\

\paragraph{\textbf{The Deligne extension}} $ $\\
For the rest of the chapter assume that all $\mct_j$ are unipotent. Define for $(z,w) \in U$ the map:
\[
\sigma_{(z,w)} : H_{\tau(z,w)} \rightarrow H_{\tau(z,w)}, \quad A \mapsto \prod_{j=1}^{r} t_j^{-N_j / (2 \pi i) } (A) \, .
\]
$\sigma_{(z,w)}$ also acts on $H^{\infty}$. Now the maps $\sigma_{(z,w)}$ map a multi-valued $\nabla$-flat section $A$ on $H$ to an elementary section on $H$ :
\[
\prod_{j=1}^{r} t_j^{-N_j /(2 \pi i)} (A) \, .
\]
It holds:
\[
\rho_{(z,w)} \circ \sigma^{-1}_{(z,w)} = \rho_{(z_1, \ldots , z_j +1 , \ldots , z_r,w)} \circ \sigma^{-1}_{(z_1, \ldots , z_j +1 , \ldots , z_r,w)} \, .
\]
Therefore we can define
\[
\chi_{\tau(z,w)} := \rho_{(z,w)} \circ \sigma^{-1}_{(z,w)} : H_{\tau(z,w)} \rightarrow H^{\infty} \, .
\]
The $\chi_{\tau(z,w)}$ map an elementary section to a constant value in $H^{\infty}$.\\
Let $\tilde{\nabla}$ be the connection for which the elementary sections are $\tilde{\nabla}$-flat:
\[
\tilde{\nabla} := \nabla + \sum_{j=1}^{r} \frac{N_j}{2 \pi i} \frac{dt_j}{t_j} \, .
\]
The elementary sections provide an extension of the vector bundle $ H \rightarrow M \setminus Z$ to a vector bundle $\tilde{H} \rightarrow M$ on $M$, the Deligne extension. The connection $\tilde{\nabla}$ gives us a trivialization of $\tilde{H} \rightarrow M$. The maps $\chi_{(t,w)}$ provide an identification of the $\tilde{\nabla}$-flat bundle $\tilde{H} \rightarrow M$ with $H^{\infty} \times M$, which we will use in the following discussion.\\

\noindent Fix $(z,w) \in U$ and consider the isomorphism $\rho_{(z,w)}: H_{\tau(z,w)} \rightarrow H^{\infty}$. Then we can shift the polarized HS on $H_{\tau(z,w)}$ to $H^{\infty}$. We denote the obtained reference polarized HS by
$(H^{\infty}, H_{\mr}^{\infty},S , F^{\bullet}_{0})$. The
space
\[
\check{D} := \{ \text{filtrations} \; F^{\bullet} \subset H^{\infty} \mid dim
F^{p} = dim F^{p}_{0}, S(F^{p},F^{w+1-p}) = 0 \}
\]
is a complex homogeneous space and a projective manifold, and the
subspace
\[
D := \{ F^{\bullet} \in \check{D} \mid F^{\bullet} \; \text{is a
part of a polarized Hodge structure} \}
\]
is an open submanifold and a real homogeneous space. It is a
classifying space for polarized Hodge structures with fixed Hodge
numbers.\\[7pt]
\noindent Now let $\mcf^{\bullet} = \bigcup_{(t,w) \in M \setminus Z} F^{\bullet}_{(t,w)}$.
The monodromies give rise to a representation of the fundamental group
\[
\phi : \pi_{1}(M \setminus Z) \longrightarrow Gl(H^{\infty}).
\]
Set $\Gamma = \phi(\pi_{1}(M\setminus Z)) = \phi(\mz^{r})$.
We have the following commutative diagram
\[
\begin{CD}
U @>\tilde{\Phi}>> D\\
@V\tau VV @VVV \\
M\setminus Z @>\Phi >> \Gamma \backslash D
\end{CD}
\]
with
\[
\tilde{\Phi}:  (z,w) \mapsto \rho_{(z,w)} (F^{\bullet}_{\tau(z,w)})\, .
\]
and where $\Phi$ is Griffiths' period mapping for the polarized VHS in question.
$ $\\
With respect to elementary sections we get the map
\begin{align}
\Psi: M \setminus Z &\longrightarrow \check{D} \notag \\
(t,w) &\mapsto \chi_{(t,w)}(F^{\bullet}_{(t,w)})\, . \notag
\end{align}
The fact that $\Psi$ takes values in $\check{D}$ rather than $D$, reflects the fact
that these sections are not real.\\[5pt]

\begin{definition}
A pair $(F^{\bullet},N)$ is said to give rise to a nilpotent orbit
if $F^{\bullet} \in \check{D}$, the endomorphism $N$ of $H_{\mr}$ is
nilpotent and an infinitesimal isometry with $N(F^{p}) \subset
F^{p-1}$, and
\[
e^{zN}F^{\bullet} \in D \quad \text{for} \;\; Im z > b\, , \quad b \in \mr .
\]
Then the set $\{ e^{zN} F^{\bullet} \mid z \in \mc \}$ is called a
nilpotent orbit.
\end{definition}
\begin{theorem}[Nilpotent Orbit Theorem]
The map $\Psi$ extends holomorphically to $\Delta^{n}$.\\ For $(w)
\in \Delta^{n-r}$, we have
\[
a(w) = \Psi(0,w) \in \check{D}\, .
\]
For any given
number $\eta$ with $0 < \eta < 1$, there exist constants
$\alpha,\beta \geq 0$, such that under the restrictions
\[
Im\, z_{i} \geq \alpha,\quad 1 \leq i \leq r, \quad \text{and} \quad
\mid w_{j}\mid  \leq \eta, \quad r+1 \leq j \leq n,
\]
the point $\exp(\sum_{i=1}^{l} z_{i} N_{i}) \circ a(w)$ lies in $D$
and satisfies the inequality
\[
d(\exp(\sum_{i=1}^{l} z_{i} N_{i}) \circ a(w), \tilde{\Phi}(z,w))
\leq (\prod_{i=1}^{r} Im(z_{i}))^{\beta} \sum_{i=1}^{l} \exp(- 2 \pi
Im(z_{i})),
\]
here $d$ denotes a $G_{\mr}$ invariant Riemannian metric on $D$.
Finally the mapping
\[
(z,w) \mapsto \exp(\sum_{i=1}^{l} z_{i} N_{i}) \circ a(w)
\]
is horizontal.
\end{theorem}
\begin{proof}
See \cite{Sch} Theorem 4.12 .
\end{proof}
$ $\\[10pt]
We denote the point $\Psi(0,0)$ as $F^{\bullet}_{lim}$.
\begin{remark}\label{3:nipo}
Observe that the pair $(F^{\bullet}_{lim}, \sum_{j=1}^{r} a_j N_j)$ gives rise to a nilpotent orbit, for any $a_j > 0$, $j= 1, \ldots , r$. Schmids Nilpotent Orbit Theorem yields
\[
\exp ( \sum_{j=1}^{r} z_{j} N_{j})\psi(0,0) \in D \quad \text{for}
\quad Im(z_{j}) > \alpha .
\]
Then
\[
\exp(z \sum_{j=1}^{r} \lambda_{j} N_{j}) \psi(0,0)  \quad
(\lambda_{j} > 0)
\]
is a nilpotent orbit, because for $Im(z) > b =
\max_{j}(\frac{\alpha}{\lambda_{j}})$ the above expression lies in
$D$.
\end{remark}
\noindent To build a PMHS on $(0,0)$ we need a nilpotent
endomorphism to define a weight filtration but the question is which
of the $N_{i}$ we should use. The following result of Cattani and
Kaplan provides an answer.
\begin{theorem}\label{3:MonoCone}
Let
\[
\mathcal{C} = \{ \sum_{j=1}^{r} \lambda_{j} N_{j} \mid \lambda_{j} >
0 \; \forall j\}.
\]
$\mathcal{C}$ is called the monodromy cone. Then each $N \in
\mathcal{C}$ defines the same monodromy weight filtration.
\end{theorem}
\begin{proof}
See, e.g. \cite{CaK} Theorem 3.3 .
\end{proof}
\noindent We need a criterion for which $N$ the tuple
$(H,H_{\mr},S,F^{\bullet},N)$ is a PMHS.
\begin{theorem}\label{3:PMHS}
$ $
\begin{enumerate}
\item The tuple $(H,H_{\mr},S,F^{\bullet},N)$ is a PMHS of weight w
$\Leftrightarrow$ the pair $(F^{\bullet},N)$ gives rise to a
nilpotent orbit.
\item If $(H,H_{\mr},S,F^{\bullet},N)$ is a PMHS with $I^{q,p} =
\overline{I^{p,q}}$, then $e^{zN}F^{\bullet} \in D$\; \text{for}\;
$Im z > 0$ .
\end{enumerate}
\end{theorem}
\begin{proof}
`$\Leftarrow$' is \cite{Sch}, Theorem 6.16. It is a consequence of the $SL_{2}$-orbit theorem. `$\Rightarrow$' is \cite{CaKSch}, Corollary 3.13. The special case (b) is proved in \cite{CaKSch}, Lemma 3.12.
\end{proof}

Now we are finally able to prove the first step in the construction
of a logarithmic Frobenius manifold out of a variation of Hodge
structures. Namely the construction of a filtration $\mcu_{\bullet}$
which is opposite to the filtration $\mcf^{\bullet}$ in vicinity of
$(0,0)$.\\
\begin{proposition}\label{3:oppFil}
Assume as above that we have given $M=\Delta^n$,
$M\setminus Z=(\Delta^*)^r\times \Delta^{n-r}$, $Z=\bigcup_{j=1}^rD_j$,
and a polarized VHS on a $\nabla$-flat bundle $H \rightarrow M\setminus Z$ with unipotent monodromies.
Let $\tilde{H} \rightarrow M$ be the Deligne extension of $H\rightarrow M\setminus Z$.
\begin{enumerate}
\item The bundle $\tilde{H} \rightarrow M$ comes equipped with a canonical
flat connection $\tilde{\nabla}$ and a $\tilde{\nabla}$-flat isomorphism
$\tilde{H} \rightarrow H^\infty\times M$.

\item The Hodge subbundles of $H\rightarrow M\setminus Z$ extend to $\tilde{H} \rightarrow M$.
The filtration $F^\bullet_{lim}$
on $\tilde{H}_0\simeq H^\infty$ is part of a PMHS.
For it one has Deligne's $I^{p,q}$. Then the filtration
$U_\bullet$ on $H^\infty$ with $U_p:=\bigoplus_{i,q:i\leq p} I^{i,q}$
is opposite to all $F^\bullet_{(t,w)}$ for $(t,w)\in M$ close to 0,
and it is $\tilde{\nabla}$-flat as well as $\nabla$-flat.
\end{enumerate}
\end{proposition}
\begin{proof}
Part (1) follows from what has been said above. For part (2), $F^{\bullet}_{lim}$ is part of a PMHS because of remark \ref{3:nipo} and theorem \ref{3:PMHS}. That $U_{\bullet}$ is opposite to $F^{\bullet}_{lim}$ is clear. We are using the $\tilde{\nabla}$-flat identification $\tilde{H} \rightarrow H^{\infty} \times M$ mentioned above to construct extensions of the $U_p$ to $\tilde{\nabla}$-flat subbundles $\calu_{p}$ on $\tilde{H}$. Because being opposite is an open property, all $F^{\bullet}_{(t,w)}$ are opposite to $U_{\bullet}$ for $(t,w)$ near 0.\\
For $N \in \calc$ we have because of lemma \ref{3:LemIpq}
\[
N: I^{p,q} \longrightarrow I^{p-1,q-1} \, .
\]
Because these N build an open cone in $\sum_{j=1}^{r} \mr \cdot N_j$ all $N_j$ are linear combinations of such $N$. Therefore we have for all $N_j$
\[
N_j : I^{p,q} \longrightarrow I^{p-1,q-1}
\]
and therefore $N_j : U_p \rightarrow U_{p-1}$. Because $U_{\bullet}$ is increasing $U_p$ is $N_j$-invariant for all $N_j$. Because we have $\tilde{\nabla} = \nabla + \sum_{j=1}^{r} \frac{N_j}{2 \pi i} \frac{dt_j}{t_j}$ we conclude that $U_{\bullet}$ is also $\nabla$-invariant.
\end{proof}

\subsection{Construction of Frobenius manifolds} $
$ \normalsize\\[15pt]
\noindent Section 3.2 is a generalization of chapter 5 of
\cite{HM}.
\begin{lemma}\label{3:LemFTSFU}
(a) The structures $(\alpha)$ and $(\beta)$ are equivalent.\\[10pt]
$(\alpha)$ A logarithmic Frobenius type structure
$((M,0),D,K,\nabla^{r},\calc,\calu,\calv,g)$ together with an
integer $w$ such that $\calu =0$ and $\calv$ is semisimple with
eigenvalues in $\frac{w}{2}+\mz$.\\[10pt]
$(\beta)$ A tuple
$((M,0),D,K,\nabla,\mcf^{\bullet},\mcu_{\bullet},w,S)$. Here $K
\rightarrow (M,0)$ is a germ of a holomorphic vector bundle;
$\nabla$ is a flat connection with logarithmic poles along $D$;
$\mcf^{\bullet}$ is a decreasing filtration by germs of holomorphic
subbundles $\mcf^{p} \subset K, p \in \mz$, which satisfies
Griffiths transversality
\[
\nabla: \calo(\mcf^{p}) \longrightarrow \omld \otimes
\calo(\mcf^{p-1});
\]
$\mcu_{\bullet}$ is an increasing filtration by subbundles
$\mcu_{p} \subset K$ which are flat outside $D$ such that
\begin{equation}\label{3:directFU}
K = \mcf^{p} \oplus \mcu_{p-1} = \bigoplus_{q} \mcf^{q} \cap
\mcu_{q};
\end{equation}
$w \in \mz$, and $S$ is a $\nabla$-flat, $(-1)^{w}$-symmetric,
nondegenerate pairing on $K$ with
\begin{align}
S(\mcf^{p},\mcf^{w+1-p}) &= 0, \label{3:SF}\\
S(\mcu_{p},\mcu_{w-1-p}) &= 0. \label{3:SU}
\end{align}
(b) One passes from $(\alpha)$ to $(\beta)$ by defining
\begin{align}
\nabla &:= \nabla^{r} + \calc, \\
\mcf^{p} &:= \bigoplus_{q \geq p} ker(\calv - (q - \frac{w}{2}) id:K
\longrightarrow K), \\
\mcu_{p} &:= \bigoplus_{q \leq p} ker(\calv - (q - \frac{w}{2}) id:K
\longrightarrow K), \\
S(a,b) &:= (-1)^{p}g(a,b) \quad \text{for} \;\; a \in \calo(\mcf^{p}
\cap \mcu_{p}), b \in \calo(K). \label{3:Sg}
\end{align}
\end{lemma}
\begin{proof}
First we prove part (b). The connection $\nabla^{r}$ and the
logarithmic Higgs field $\calc$ are maps $\calo(K) \rightarrow \omld
\otimes \calo(K)$. Then the flatness $\nabla^{r}$ means
$(\nabla^{r})^{2}=0$, the Higgs field $\calc$ satisfies $\calc^{2}
=0$, and the potentiality condition means $\nabla^{r}(\calc) =
\nabla^{r} \circ \calc + \calc \circ  \nabla^{r} =0$. Therefore
$\nabla^{2} = (\nabla^{r} + \calc)^{2} =0$, the connection $\nabla$
is flat and has logarithmic poles along $D$.\\
The filtrations $\mcf^{\bullet}$ and $\mcu_{\bullet}$ obviously
satisfy (\ref{3:directFU}) and
\[
\mcf^{p} \cap \mcu_{p} = ker(\calv - (p- \frac{w}{2}) id:\; K
\longrightarrow K).
\]
The connection $\nabla^{r}$ maps $\calo(\mcf^{p} \cap \mcu_{p})$ to
itself because $\calv$ is $\nabla^{r}$-flat. Because $\calu = 0$ we
have $ [ \calc ,\calv] = \calc$, which is equivalent to $\calc_{X},
X \in \dmldc$ mapping $\calo(\mcf^{p} \cap \mcu_{p})$  to
$\calo(\mcf^{p-1} \cap \mcu_{p-1})$. Therefore $\mcu_{\bullet}$ is
$\nabla$-flat and
$\mcf^{\bullet}$ satisfies Griffiths transversality. Obviously $\nabla^r$ and $\calc$ have logarithmic poles along $D$.\\[10pt]
The condition (\ref{1:gV}) says that for $ a \in \calo(\mcf^{p} \cap
\mcu_{p}), b \in \calo(\mcf^{q} \cap \mcu_{q})$
\[
S(a,b) = (-1)^{p} g(a,b) = 0 \quad \text{if} \;\; p+q \neq w\, .
\]
Therefore $S$ is $(-1)^{w}$-symmetric and satisfies (\ref{3:SF}) and
(\ref{3:SU}).

\noindent To show that it is $\nabla$-flat we have to choose a representant of $((M,0),D,K,\nabla^r)$
which we denote by the same letters. Now let $U \subset M\setminus D$ open and simply connected.
For $X \in \dmldc(U)$ and
$\nabla^{r}$-flat sections $a \in \calo(\mcf^{p} \cap \mcu_{p})(U),\, b
\in \calo(\mcf^{q} \cap \mcu_{q})(U)$ we have
\begin{align}
(\nabla_{X}S)(a,b) &= X\,S(a,b) - S(\nabla^{r}_{X}a + \calc_{X} a,b)
- S(a,\nabla^{r}_{X}b + \calc_{X} b) \notag \\
&= 0 - (-1)^{p+1} g(\calc_{X}a,b) - (-1)^{p}g(a,\calc_{X}b)=0.
\notag
\end{align}
This shows (b).\\[10pt]
To pass from $(\beta)$ to $(\alpha)$ we notice that $\calv$ is
uniquely determined by
\[
\mcf^{p} \cap \mcu_{p} = ker(\calv -(p-\frac{w}{2})id: K
\longrightarrow K)
\]
and $g(a,b), a \in \calo(\mcf^{p} \cap \mcu_{p}), b \in \calo(K)$ is
defined by (\ref{3:Sg}).\\
We can decompose $\nabla$ into $\nabla^{r}$ mapping $\calo(\mcf^{p}
\cap \mcu_{p})$ to itself and $\calc_{X}$ for $X \in \dmldc$ mapping
$\calo(\mcf^{p} \cap \mcu_{p})$ to $\calo(\mcf^{p-1} \cap
\mcu_{p-1})$ respectively. This is possible because
$\mcu_{\bullet}$ is flat and $\mcf^{\bullet}$ satisfies Griffiths
transversality.\\[10pt]
Now we decompose $\nabla^{2}$ into its components. Thus we have
$(\nabla^{r})^{2} = 0$, $\nabla^{r} \circ \calc + \calc \circ
\nabla^{r} = 0$ and $\calc^{2} = 0$. The second is equivalent to
condition (\ref{D :FTS1}) in the definition of logarithmic Frobenius
structure and the
third shows that $\calc$ is a Higgs field with logarithmic poles.\\
Condition (\ref{D :FTS2}) is satisfied because $\calu = 0$ and $[
\calc , \calv
] = \calc$ per construction of $\calc$.\\
The condition $g(\calc_{X} a,b) = g(a, \calc_{X}b)$ follows from the
very same calculation above.\\
$g(\calv a,b) = - g(a, \calv b)$ follows right from the definitions.
\end{proof}

\begin{remark}
If one adds  in $(\beta)$ a real structure with suitable conditions
one obtains a germ of a variation of polarized Hodge structures of
weight $w$.
\end{remark}
\begin{definition}\label{3:H2Frob}
$ $\\
(a) A germ of an $H^{2}$-generated variation of filtrations of
weight $w$, which has logarithmic poles at $D$, is a tuple
$((M,0),D,K,\nabla, \mcf^{\bullet},w)$ with the following
properties: $w$ is an integer, $K \rightarrow (M,0)$ is a germ of a
holomorphic vector bundle with flat connection $\nabla$ with
logarithmic poles along $D$ and a variation of filtrations
$\mcf^{\bullet}$ which satisfies Griffiths transversality and
\begin{align}
& 0 = \mcf^{w+1} \subset \mcf^{w} \subset \ldots \subset K\, , \notag \\
& rk \mcf^{w} =1 , \; rk \mcf^{w-1} = 1 + dim M \geq 2\, . \notag
\end{align}
Griffiths transversality and flatness of $\nabla$ give a Higgs field
$\calc$ on the graded bundle $\bigoplus_{p} \mcf^{p} / \mcf^{p+1}$
with commuting endomorphisms
\[
\calc_{X} = [\nabla_{X}]: \calo(\mcf^{p} / \mcf^{p+1})
\longrightarrow \calo(\mcf^{p-1}/ \mcf^{p}) \quad \text{for} \; X
\in \dmldc\, .
\]
$ $\\
$\mathbf{H^{2}}$\textbf{-generation condition:} The whole module
$\bigoplus_{p}\calo(\mcf^{p}/\mcf^{p+1})$ is generated by
$\calo(\mcf^{w})$ and its images under iterations of the maps
$\calc_{X}, X \in \dmldc$.\\[15pt]
(b) A pairing and an opposite filtration for an $H^{2}$-generated
variation of filtrations $((M,0),D,$ $K,\nabla, \mcf^{\bullet})$ of
weight $w$ are a pairing $S$ and a filtration $\mcu_{\bullet}$ as in
lemma (\ref{3:LemFTSFU}).
\end{definition}
$ $\\
\begin{definition}
An $H^{2}$-generated germ of a logarithmic Frobenius manifold of
weight $w \in \mathbb{N}_{\geq 1}$ is a germ $((M,0),D,\circ,e,E,g)$
of a logarithmic Frobenius manifold with the properties (I) and (II)
below and with
\[
E\mid_{t=0} = 0.
\]
Let $\nabla^{g}$ be the Levi-Civita connection of $((M,0),g)$. It has logarithmic poles along $D$ by proposition \ref{1:LClp}. The
endomorphism
\begin{align}
\nabla^{g}E: \dmldc_{0} &\rightarrow \dmldc_{0}\, , \notag \\
X &\mapsto \nabla^{g}_{X}E \notag
\end{align}
acts on the space of $\nabla^{g}$-flat logarithmic vector fields.\\
In particular $e \in ker(\nabla^{g}E -id)$.\\
(I) It acts semisimply with eigenvalues $\lbrace 1,0, \ldots ,
-(w-1) \rbrace$.\\
It turns out that then the multiplication on the algebra
$\dmld_{0}$ respects the grading
\[
\dmld_{0} = \bigoplus_{p=0}^{w} ker(\nabla^{g}E - (1-p) id :
\dmld_{0} \longrightarrow \dmld_{0})
\]
and that $Lie_{E}(g) = (2-w) \cdot g$.\\
(II) $H^{2}$-generation condition: The algebra $\dmld_{0}$ is
generated by \\
$ker(\nabla^{g} E : \dmld_{0} \longrightarrow \dmld_{0})$.
\end{definition}
$ $\\
\begin{theorem}\label{3:corespondFMH2}
There is a one-to-one correspondence between the structures in
$(\alpha),(\beta)$ and $(\gamma)$.\\[10pt]
$(\alpha)$ A logarithmic Frobenius type structure
$((M,0),D,K,\nabla^{r},\calc,\calu,\calv,g)$ with $\calu = 0$ and
with a fixed section $\xi \in \calo(K)_{0}$ which is $\nabla^r$-flat on $M \setminus D$ and which
satisfies the conditions (IC),(GC) and (EC) in theorem
(\ref{U:General}).\\[10pt]
$(\beta)$ A germ of an $H^{2}$-generated variation of filtrations
with logarithmic pole along $D$
$((M,0),$ $D, K, \nabla,\mcf^{\bullet}, w, S, \mcu_{\bullet})$ of weight $w
\in \mathbb{N}_{\geq 1}$ with pairing and opposite filtration and
with a fixed generator $\xi \in \calo(\mcf^{w})$ such that $\nabla \xi \in \calo(\calu_{w-1}) \otimes \dmldc$.\\[10pt]
$(\gamma)$ An $H^{2}$-generated germ of a logarithmic Frobenius
manifold $((\tilde{M},0),\tilde{D},\circ,e,$ $E,\tilde{g})$ of weight
$w \in
\mathbb{N}_{\geq 1}$.\\[10pt]
One passes from $(\alpha)$ to $(\beta)$ by lemma (\ref{3:LemFTSFU})
and from $(\alpha)$ to $(\gamma)$ by theorem (\ref{U:General}). One
passes from $(\gamma)$ to $(\alpha)$ by defining
\begin{align}
M &:= \lbrace t \in \tilde{M} \mid E\mid_{t} = 0 \rbrace, \notag \\
D &:= \tilde{D} \cap (M,0)\, , \notag \\
K &:= Der_{\tilde{M}}(log \tilde{D})\mid_{(M,0)} \notag
\end{align}
with the canonical logarithmic Frobenius type structure, and $\xi :=
e\mid_{(M,0)}$. The eigenvector condition $(EC)$ is $\calv \xi =
\frac{w}{2} \xi$.
\end{theorem}
\noindent In contrast to the corresponding theorem in \cite{HM} where they could simply extend $\xi\mid_0$ $\nabla^r$-flat, we have to demand the existence of a section $\xi$ which is flat on $M \setminus D$.
\begin{proof}
$(\alpha) \Rightarrow (\beta)$.\\
We have to show that the conditions in lemma \ref{3:LemFTSFU} $(\alpha)$ are satisfied.
For that we have to show that for some $w \in \mathbb{N}_{\geq 1}$ the
endomorphism $\calv +\frac{w}{2}id$ is semisimple with eigenvalues
in $\lbrace 0,1, \ldots , w \rbrace$. For the conditions on $\mcf^{\bullet}$ in definition \ref{3:H2Frob} we have to show $\calv \xi =
\frac{w}{2} \xi$.\\
The eigenvector condition says $\calv \xi = \frac{d}{2} \xi$ for
some $d \in \mc$. Condition (\ref{D :FTS2}) reads as $[\calc ,
\calv] = \calc$. This together with the generation condition (GC)
shows that $\calv$ acts semisimply on $K_{0}$ with eigenvalues in
$\frac{d}{2} + \mz_{\leq 0}$ and that $ker(\calv - \frac{d}{2}id) =
\mc \cdot \xi$. The injectivity condition (IC) implies $dim\,
ker(\calv - (\frac{d}{2} -1)id) = dim M > 0$. Condition (\ref{1:gV})
says that the eigenvalues of $\calv$ are in $(\frac{d}{2} +
\mz_{\leq 0}) \cap -(\frac{d}{2} + \mz_{\leq 0})$. Therefore $ d \in
\mathbb{N}$. Define $w := d$. The conditions (IC) and (GC) show that
one passes from $(\alpha)$ to
$(\beta)$ by lemma (\ref{3:LemFTSFU}). $\nabla \xi\mid_{M\setminus D} \in \calo(\calu_{w-1})$ follows from $\nabla^r \xi \mid_{M \setminus D} = 0$.\\

\noindent $(\beta) \Rightarrow (\alpha)$\\
Lemma (\ref{3:LemFTSFU}) gives a Frobenius type structure. It holds $\nabla^r \xi\hspace{-4pt}\mid_{M\setminus D} = 0$ because $\nabla \xi \in \calo(\calu_{w-1}) \otimes \dmldc$. We have
to show (IC),(GC) and (EC).\\
(IC) follows from the condition $rk \mcf^{w-1} = 1 +dim M$ and from the $H^2$-generation condition. $(GC)$
follows from the $H^{2}$-generation condition. Finally (EC) follows
from the characterization of $\calv$ as $\mcf^{p} \cap \mcu_{p} =
ker (\calv - (p - \frac{w}{2})id :K \longrightarrow K)$,$\;\;rk
\mcf^{w} = 1$ and $K = \bigoplus_{q} \mcf^{q} \cap
\mcu_{q}$.\\[5pt]
\noindent $(\alpha) \Rightarrow (\gamma)$\\
Theorem (\ref{U:General}) applied to $(\alpha)$ gives a logarithmic
Frobenius manifold $((\tilde{M},0),D,$ $\circ,e,E,\tilde{g})$, an
embedding $i: (M,0) \rightarrow (\tilde{M},0)$ and an
isomorphism of logarithmic Frobenius type structures $j : K \rightarrow Der_{\tilde{M}}(log
\tilde{D})\mid_{i(M)}$. It
maps $\calv$ to the restriction of $\nabla^{g}E- \frac{2-d}{2} id$
on $Der_{\tilde{M}}(log \tilde{D})\mid_{i(M)}$. Therefore one
obtains an $H^{2}$-generated germ of a Frobenius manifold of weight
$w$. In suitable flat coordinates the Euler field of this Frobenius
manifold takes the form
\[
E = \sum_{i=1}^{dim \tilde{M}} d_{i} t_{i} \frac{\p}{\p t_{i}}
\]
with $d_{i} \in \lbrace 1,0, \ldots , -(w-1) \rbrace$ and
\[
\sharp ( i \mid d_{i} =0 ) = dim\; ker(\nabla^{g} E) = dim\; ker(\calv
-(\frac{d}{2} -1)id) = dim M.
\]
Observe that for logarithmic coordinates $d_{i}=0$, because $\calu= 0$ on $M$.\\
Therefore $i(M)
= \lbrace t \in \tilde{M} \mid E\mid_{t} =0 \rbrace$.\\[2pt]
$(\gamma) \Rightarrow (\alpha)$\\
One passes from $(\gamma)$ to $(\alpha)$ as described above. (IC) is
fulfilled because $\xi = e$. (GC) is because the Frobenius manifold
is $H^{2}$-generated. (EC) is also from the definition.
\end{proof}

\subsection{Application to VPHS}$ $ \normalsize$ $\\[3pt]
In this section we apply the machinery developed above to VPHS and investigate
some properties of the resulting Frobenius type structures.\\[-12pt]
\begin{proposition}\label{3:VPHSFrob}
Assume as above that we have given $M= \Delta^{n}$, $M \setminus Z = (\Delta^\ast)^r \times \Delta^{n-r}$,
$Z= \bigcup_{j=1}^{r} D_j$, and a VPHS on a $\nabla$-flat bundle $H \rightarrow M \setminus Z$ with
unipotent monodromies. Let $\tilde{H} \rightarrow M$ be the Deligne extension of
$H \rightarrow M \setminus Z$.\\
Then one can construct a logarithmic Frobenius type structure and a $(logD-trTLEP(w))$-structure out of
this data.
\end{proposition}
\begin{proof}
Set $K:= \tilde{H}$. Construct an increasing filtration by $\nabla$-flat subbundles $\calu_p$
with the help of lemma \ref{3:oppFil}, which is opposite to the filtration $\mcf_{\bullet}$.
Now lemma \ref{3:LemIpq} gives us $S(I^{p,q},I^{r,s}) = 0$ for $(r,s) \neq (w-p,w-q)$ and
because $U_p = \sum_{i,q:i \leq p} I^{i,q}$ we see that condition (\ref{3:SU})
$(S(\calu_p,\calu_{w-1-p})= 0)$ is fulfilled. Now we have data as in lemma
$\ref{3:LemFTSFU}(a)(\beta)$ and therefore a logarithmic Frobenius type structure
$((M,0),Z,K,\nabla^{r},\calc,$ $\calu,\calv,g,w)$ with $\calu =0$ and
$\calv$ semisimple with eigenvalues in $\frac{w}{2} + \mz$.
Lemma \ref{FTSlogDoneone} provides a $(logD-trTLEP(w))$-structure.
\end{proof}

\begin{proposition}\label{3:D2}$ $\\[-9pt]
\begin{enumerate}
\item The connection $\nabla^r$ has trivial monodromy.
\item It can be extended holomorphically from $H \rightarrow M \setminus Z$ to $K \rightarrow M$.
\end{enumerate}
\end{proposition}
\begin{proof}
Recall the proof of lemma \ref{3:LemFTSFU} part $(\beta) \rightarrow (\alpha)$. There the map
\[
\nabla: \calo(\mcf^p\cap \calu_p) \rightarrow (\calo(\mcf^p \cap \mcu_p) \oplus \calo(\mcf^{p-1} \cap \calu_{p-1})) \otimes  \dmldc ,
\]
was decomposed into $\nabla^r$ mapping to $\calo(\mcf^p \cap \mcu_p) \otimes  \dmldc$ and $\calc$ mapping to
$\calo(\mcf^{p-1} \cap \calu_{p-1}) \otimes  \dmldc$ respectively. A different characterization of
$\nabla^r$ and $\calc$ would be the following: All $\calu_p$ are $\nabla$-flat subbundles of $K$. Now
$\nabla$ induces flat connections $\nabla^{(U)}$ on the quotients $\calu_p / \calu_{p-1}$. Because
$\mcf^{\bullet}$ and $\calu_{\bullet}$ are opposite we have canonical isomorphisms
$\mcf^p \cap \calu_p \simeq \calu_p / \calu_{p-1}$ and therefore
\[
\bigoplus_p \calu_p / \calu_{p-1} \simeq \bigoplus_p \mcf^{p} \cap \calu_p = K\, .
\]
Under this isomorphism $\nabla^{(U)}$ goes over to $\nabla^r$ and $\calc:= \nabla -\nabla^r$.
Because of $\mct_j = e^{N_j}$ and $N_j: \calu_p \rightarrow \calu_{p-1}$, $\mct_j$ induces the identity on $\calu_p / \calu_{p-1}$. Therefore also $\nabla^r$ has trivial monodromy. This shows (1).\\
For (2) recall that $\calu_p$ is the Deligne extension of its restriction $\calu_p \mid_{M\setminus Z}$ on
$M \setminus Z$ with respect to $\nabla$. Therefore $\calu_p / \calu_{p-1}$ is the Deligne extension
with respect to $\nabla^{(U)}$ of its restriction to $M \setminus Z$. Because of that and (1),
$\nabla^{(U)}$ has an extension to $M$. This shows (2).
\end{proof}
We will now give conditions on the PMHS on $\tilde{H}_0$ so that the logarithmic Frobenius type structure constructed in proposition \ref{3:VPHSFrob} gives rise to a logarithmic Frobenius manifold.
\begin{theorem}\label{3:VPHSFrob2}
Assume that we have given $M = \Delta^n $ and $ M \setminus Z = (\Delta^{\ast})^n$ and a VPHS on a $\nabla$-flat bundle $H \rightarrow M \setminus Z$ with unipotent monodromies $\mathcal{T}_i$. Let $\tilde{H} \rightarrow M$ be the Deligne extension and $F^{\bullet}_{lim}$ be part of the limiting MHS on $\tilde{H}_0$ with
\begin{align}
&0 = \,F^{w+1}_{lim} \subset F^w_{lim} \subset \ldots \subset K \notag \\
\text{and}\quad &dim\, F^w_{lim} =1 , \;\; dim\, F^{w-1}_{lim} = 1 + dim M . \notag
\end{align}
Let $\hspace{-2pt}N_i\hspace{-2pt}$ be the logarithm of the unipotent part of $\mathcal{T}_i\hspace{-2pt}$ acting on $\hspace{-1pt}\tilde{H}_0$. If the vector space $\bigoplus_p\hspace{-2pt} F^p_{lim}/\hspace{-1pt}F^{p+1}_{lim}$ is generated by $F^w_{lim}$ and its images under iterations of the linear maps $N_i$, then we can construct an $H^2$-generated germ of a logarithmic Frobenius manifold out of these data.
\end{theorem}
\begin{proof}
Pick a nonzero vector $v \in F^{w}_{lim}$. Recall proposition \ref{3:VPHSFrob} where we constructed a logarithmic Frobenius type structure $((M,0),Z,K, \nabla^r,\calc,\calu,\calv,g,w)$ out of the VPHS. Because of proposition \ref{3:D2} we can extend this vector $v$ to a $\nabla^r$-flat section $\xi \in \calo(K)$. The conditions on $F^{\bullet}_{lim}$ and $v \in F^w_{lim}$ give the conditions $(IC)$ and $(GC)$ in theorem \ref{U:General} resp. theorem \ref{3:corespondFMH2} with respect to the section $\xi \in \calo(K)$. The condition (EC) follows from the construction of the logarithmic Frobenius type structure in lemma \ref{3:LemFTSFU} and the existence of an opposite filtration $U_{\bullet}$ on $\tilde{H}_0 = K_0$ (cf. \ref{3:oppFil}).
\end{proof}
\begin{remark}
From the generation conditions and the properties of $F^{\bullet}_{lim}$ follows that the PMHS on $\tilde{H}_0$ is Hodge-Tate and split over $\mathbb{R}$. This can be easily seen by choosing $v$ to to be real and by noticing that every $N_j$ is a $(-1,-1)$-morphism defined on $\tilde{H}_{\mathbb{R},0}$.
\end{remark}

\end{document}